\documentclass{lmcs}
\pdfoutput=1

\usepackage{lastpage}
\lmcsdoi{15}{1}{28}
\lmcsheading{}{\pageref{LastPage}}{}{}%
{Dec.~11,~2018}{Mar.~21,~2019}{}

\usepackage[utf8]{inputenc}
\usepackage{csquotes}
\usepackage[english]{babel}  
\usepackage[T1]{fontenc} %

\usepackage{xspace}
\usepackage{relsize}

\usepackage{amsmath}
\usepackage{amsfonts}
\usepackage{amssymb}
\usepackage{amsthm}

\usepackage[all,2cell]{xy}\SelectTips{cm}{10}\SilentMatrices\UseAllTwocells 
\usepackage[english]{babel}
\usepackage{xcolor}
\definecolor{dkblue}{rgb}{0,0.1,0.5}
\definecolor{lightblue}{rgb}{0,0.5,0.5}
\definecolor{dkgreen}{rgb}{0,0.4,0}
\definecolor{dk2green}{rgb}{0.4,0,0}
\definecolor{dkviolet}{rgb}{0.6,0,0.8}
\definecolor{DarkBlue}{rgb}{0,0.08,0.45}

\usepackage{xspace}
\usepackage{verbatim}
\usepackage{listings}

\newcommand{\fixedtextsf}[1]{\ensuremath{\mathsf{#1}}}
\theoremstyle{definition}
\newtheorem{form}[thm]{Code}

\DeclareMathOperator{\true}{true}
\DeclareMathOperator{\false}{false}

\DeclareMathOperator{\Id}{Id}

\DeclareMathOperator{\id}{id}

\DeclareMathOperator{\lift}{lift}

\DeclareMathOperator{\App}{App}
\DeclareMathOperator{\Abs}{Abs}
\DeclareMathOperator{\Var}{Var}
\DeclareMathOperator{\app}{app}
\DeclareMathOperator{\abs}{abs}

\DeclareMathOperator{\Rep}{Model}

\DeclareMathOperator{\rmlift}{mlift}

\DeclareMathOperator{\Succ}{succ}

\DeclareMathOperator{\zero}{zero}

\DeclareMathOperator{\ULCop}{\mathsf{ULC}}
\DeclareMathOperator{\substop}{subst}

\DeclareMathOperator{\LModop}{\mathsf{LMod}}
\DeclareMathOperator{\LRModop}{\mathsf{LRMod}}
\DeclareMathOperator{\RModop}{\mathsf{RMod}}

\DeclareMathOperator{\dom}{dom}
\DeclareMathOperator{\cod}{cod}

\DeclareMathOperator{\SLC}{\mathsf{TLC}}

\DeclareMathOperator{\SUCC}{Succ}

\newcommand{\LMod}[3]{\LModop_{#1}({#2},{#3})}
\newcommand{\LRMod}[3]{\LRModop_{#1}({#2},{#3})}
\newcommand{\RMod}[2]{\RModop({#1},{#2})}

\newcommand{\Set}{\ensuremath{\mathsf{Set}}}
\newcommand{\SET}{\Set}
\newcommand{\PO}{\ensuremath{\mathsf{Pre}}}
\newcommand{\Fin}{\ensuremath{\mathsf{Fin}}}

\newcommand{\nt}{\ensuremath{\Longrightarrow}}

\newcommand{\SigMon}[1]{{#1}\text{-}\mathsf{Mnd}}
\newcommand{\SigRMon}[1]{{#1}\text{-}\mathsf{RMnd}}

\newcommand{\kl}[2][]{\sigma^{#1}\left({#2}\right)}
\newcommand{\mkl}[2][]{\varsigma^{#1}({#2})}

\newcommand{\we}{\eta}
\newcommand{\comp}[2]{\ensuremath{{#2} \circ {#1}}}

\newcommand{\T}{\mathcal{T}}
\newcommand{\NN}{\mathbb{N}}
\newcommand{\C}{{\mathcal{C}}}
\newcommand{\D}{\mathcal{D}}
\newcommand{\E}{\mathcal{E}}
\newcommand{\X}{\mathcal{X}}

\newcommand{\family}[2]{\ensuremath{{#1}^{#2}}}

\newcommand{\TS}[1]{\family{\Set}{#1}}
\newcommand{\TP}[1]{\family{\PO}{#1}}
\newcommand{\TDelta}[1]{\Delta^{#1}}

\newcommand{\shift}[2]{{}_{#2}{#1}}

\newcommand{\subst}{\substop}

\newcommand{\LC}{\ULC}
\newcommand{\ULC}{\ensuremath{\ULCop}}

\newcommand{\TLC}{\SLC}
\newcommand{\TLCB}{\ensuremath{\SLC_{\mathsf{\beta}}}}

\newcommand{\PCF}{\ensuremath{\mathsf{PCF}}\xspace}

\newcommand{\True}{\ensuremath{\mathbf{T}}}
\newcommand{\False}{\ensuremath{\mathbf{F}}}

\newcommand{\Bool}{o}
\newcommand{\Nat}{\iota}
\newcommand{\zeroqu}{\ensuremath{\zero ?}}
\newcommand{\PCFFix}{\ensuremath{\mathbf{Fix}}}
\newcommand{\PCFSucc}{\ensuremath{\mathbf{Succ}}}
\newcommand{\PCFn}[1]{\ensuremath{\mathbf{#1}}}

\newcommand{\PCFZerotest}{\ensuremath{\mathbf{Zero?}}}
\newcommand{\PCFcond}[1]{\ensuremath{\mathbf{cond_{#1}}}}
\newcommand{\PCFTrue}{\ensuremath{\mathbf{T}}}
\newcommand{\PCFFalse}{\ensuremath{\mathbf{F}}}
\newcommand{\PCFPred}{\ensuremath{\mathbf{Pred}}}

\newcommand{\PCFar}{\Rightarrow}
\newcommand{\TLCar}{\leadsto}
\newcommand{\TLCTYPE}{\ensuremath{T_{\SLC}}}

\newcommand{\bind}[2]{{#1}\gg\hspace{-.5ex}={#2}}
\newcommand{\coq}{\fixedtextsf{Coq}\xspace}

\newcommand{\fibre}[2]{[{#1}]_{#2}}

\newcommand{\retyping}[1]{\ensuremath{\vec{#1}}}
\newcommand{\vectorletter}[1]{\mathbf{#1}}
\newcommand{\init}[1]{\ensuremath{\hat{#1}}}

\newcommand{\inr}{\mathrm{inr}}
\newcommand{\inl}{\mathrm{inl}}

\newcommand{\sourceurl}{\url{https://github.com/benediktahrens/monads}}

\begin{document}

\title{Initial Semantics for Reduction Rules}
\author{Benedikt Ahrens}
\address{School of Computer Science, University of Birmingham, United Kingdom}
\email{b.ahrens@cs.bham.ac.uk}

\begin{abstract}

 We give an algebraic characterization of the syntax and operational semantics 
of a class of simply--typed languages, such as the language $\PCF$: 
we characterize simply--typed syntax with variable binding and equipped with reduction rules
via a universal property, namely as the initial object of some category of models.
For this purpose, we employ techniques developed in two previous works: in the first work we
model syntactic translations between languages \emph{over different sets of types} as initial morphisms
in a category of models. In the second work we characterize untyped syntax \emph{with reduction rules}
as initial object in a category of models.
In the present work, we combine the techniques used earlier in order to
 characterize
simply--typed syntax with reduction rules as initial object in a category. 
The universal property yields an operator which allows to specify
translations---that are semantically faithful by construction---%
between languages over possibly different sets of types.

As an example, we upgrade a translation from $\PCF$ to the untyped lambda calculus, 
given in previous work, to account for reduction in the source and target.
Specifically, we specify a reduction semantics in the source and target language through suitable
rules. By equipping the untyped lambda calculus with the structure of a model of $\PCF$, 
initiality yields a translation from $\PCF$ to the lambda calculus, that is faithful with respect to
the reduction semantics specified by the rules.

This paper is an extended version of an article published in the proceedings of WoLLIC 2012.

\end{abstract}

\maketitle

\section{Introduction}

The syntactic structure of languages is captured algebraically by the notion of \emph{initial algebra}:
Birkhoff \cite{birkhoff1935} first characterizes \emph{first-order} syntax with equalities between terms as initial object in some category
(without mention of the word \enquote{category}, of course).
Algebraic characterizations of \emph{higher--order} syntax, i.e., of syntax with binding of variables, 
were given  by several people independently in the late 1990's, using different
approaches to the representation of binding.
Some of these approaches were generalized later to deal with \emph{simple typing} as well as with \emph{equations} between terms, such as 
$\beta$- and $\eta$-conversion for functions.

In this work we give an algebraic characterization of simply--typed languages equipped with \emph{directed} equations.
Directed equations provide a more faithful model for computation---or \emph{reduction}---than equations and 
emphasize the \emph{dynamic} point of view of reductions rather than the static one given by equations,
as pointed out in \cite{proofs_and_types}.

The languages we consider in this work are \emph{purely functional}; 
in particular, we do not treat memory, environment, effects or exceptions.
The rewrite rules we consider are \emph{unconditional} rewrite rules, such as $\beta$- and $\eta$-reduction.
We do not treat conditional rewrite rules as needed, e.g., in the definition of bisimulation 
for process calculi.

By including reduction rules into our treatment of the syntax of programming languages, 
we obtain a mechanism---a \emph{recursion operator}---that allows us to define translations between languages
that are \emph{semantically faithful} by construction.
A semantically faithful translation $f$ is a translation that commutes with reduction: 
it maps two terms $t$ and $t'$, where $t$ reduces to $t'$, to terms
$f(t)$ and $f(t')$ such that $f(t)$ reduces to $f(t')$ in the target language---an important property 
of any \emph{reasonable} translation between programming languages.

Our definitions ensure furthermore that translations specified 
via this operator are by construction compatible with substitution (through our use of (relative) monads) and typing, in addition to reduction.

As an example, we specify a translation from $\PCF$ to the untyped lambda calculus $\ULC$ using this
category--theoretic iteration operator. This translation is by construction faithful with respect to
reduction in $\PCF$ and $\ULC$.
This example is verified formally in the proof assistant \coq \cite{coq}.
The \coq files as well as documentation are available online\footnote{\sourceurl}.

\subsection{Difference to the conference version of this article}

The present work is an extended version of a previously published work \cite{ahrens_wollic2012}. 
In that previous work, the main theorem \cite[Thm.\ 44]{ahrens_wollic2012} is stated, but no proof is given.
In the present work, we review the definitions given in the earlier work
and present a proof of the main theorem.
Afterwards, we explain in detail the formal verification in the proof assistant \coq
\cite{coq} of an instance of this theorem, for the simply--typed programming language $\PCF$.
Finally, we illustrate the iteration operator coming from initiality by specifying 
 an executable certified translation in \coq from $\PCF$ to the untyped lambda calculus.

\subsection{Summary}

We define a notion of \emph{2--signature} which allows the specification of
the \emph{types and terms} of a programming language, 
as well as its \emph{operational semantics} in form of reduction rules for terms.
Specifically, such a 2--signature consists, first, of a 1--signature, say, $(S,\Sigma)$,
where $S$ is a signature for types, and $\Sigma$ is a signature for terms over the types specified by $S$.
Second, on such a 1--signature $(S,\Sigma)$, reductions rules are specified through 
a set $A$ of \emph{directed equations}---also called \emph{reduction rules} or \emph{rules} for short---over $(S,\Sigma)$.
To any 1--signature $(S,\Sigma)$ we associate a category of \emph{models of $(S,\Sigma)$}. 
Given a 2--signature $((S,\Sigma),A)$, the rules of $A$ give rise to a \emph{satisfaction}
predicate on the models of $(S,\Sigma)$, and thus specify the full subcategory of models of $(S,\Sigma)$
that satisfy the rules of $A$.
We call this subcategory the \emph{category of models of $((S,\Sigma),A)$}.
Our main theorem states that this category has an initial object---%
the programming language associated to $((S,\Sigma),A)$---,
which integrates the types and terms generated by $(S,\Sigma)$, equipped with the reduction relation
generated by the rules of $A$.

This characterization of syntax with reduction rules is given in two steps:
\begin{enumerate}
 \item 
At first \emph{pure syntax} is characterized as initial object in some category.
Here we use the term \enquote{pure} to express the fact that no semantic aspects such as reductions on 
terms are considered. As will be explained in Section\ \ref{sec:intro_pure_syntax}, this characterization is actually a consequence of an earlier result 
 \cite{ahrens_ext_init}.
 \item  Afterwards we consider \emph{directed equations} specifying \emph{reduction rules}. Given a family of reduction rules
          for terms, we build up on the preceding result to give an algebraic characterization of 
          syntax \emph{with reduction}.
        Directed equations for \emph{untyped} syntax are considered in earlier work \cite{ahrens_relmonads}; 
  in the present work, the main result of that earlier work is carried over to \emph{simply--typed} syntax.
\end{enumerate}

\noindent
In summary, the merit of this work is to give an algebraic characterization of simply--typed syntax \emph{with reduction rules},
building up on such a characterization for \emph{pure} syntax given earlier \cite{ahrens_ext_init}.
Our approach is based on relative monads
as defined in \cite{DBLP:journals/corr/AltenkirchCU14}  
from the category $\Set$ of sets to the category $\PO$ of preordered sets. 
Compared to traditional monads, relative monads allow for different categories as domain and codomain.

We now explain the above two points in more detail:

\subsubsection{Pure Syntax}\label{sec:intro_pure_syntax}

A 1--signature $(S,\Sigma)$ is a pair which specifies the types and terms of a language, respectively. Furthermore, it associates 
a type to any term.
To any 1--signature $(S,\Sigma)$ we associate a category $\Rep^{\Delta}(S,\Sigma)$ of \emph{models} of $(S,\Sigma)$,
where a model of $(S,\Sigma)$ is a pair of a model $T$ of the types specified by $S$ and a model of $\Sigma$ over $T$.
Models of $\Sigma$ are relative monads, from (families of) sets to (families of) preordered sets, equipped with some 
extra structure.

This category has an initial object (cf.\ Lemma\ \ref{lem:init_no_eqs_typed}), which integrates the types and terms
freely generated by $(S,\Sigma)$. We call this object the \emph{(pure) syntax associated to $(S,\Sigma)$}. 
As mentioned above, we use the term \enquote{pure} to distinguish this initial object from 
the initial object associated to a \emph{2--signature}, 
which gives an analogous characterization of syntax \emph{with reduction rules} (cf.\ below).

Initiality for pure syntax is actually a consequence of a related initiality theorem proved in another work \cite{ahrens_ext_init}:
in that work, we associate, to any signature $(S,\Sigma)$, a category $\Rep(S,\Sigma)$ of models of $(S,\Sigma)$, where a model is built from 
a traditional monad on (families of) sets instead of a relative monad as above.
We connect the corresponding categories by exhibiting a pair of adjoint functors (cf.\ Lemma\ \ref{lem:init_no_eqs_typed}) 
between the category $\Rep^{\Delta}(S,\Sigma)$ of models of $(S,\Sigma)$ and the category $\Rep(S,\Sigma)$,
  
\begin{equation*}   
       \begin{xy}
        \xymatrix@C=4pc{
                  **[l] \Rep(S,\Sigma) \rtwocell<5>_{U_*}^{\Delta_*}{'\bot} & **[r] \Rep^{\Delta}(S,\Sigma)
}
       \end{xy} \enspace .
\end{equation*}
We thus obtain an initial object in our category $\Rep^{\Delta}(S,\Sigma)$ using the fact that left adjoints are co\-con\-tin\-uous:
the image under the functor 
$\Delta_* : \Rep(S,\Sigma) \to \Rep^{\Delta}(S,\Sigma)$ 
of the initial object in the category $\Rep(S,\Sigma)$ is initial in $\Rep^{\Delta}(S,\Sigma)$.

\subsubsection{Syntax with Reduction Rules}

Given a 1--signature $(S,\Sigma)$, an \emph{$(S,\Sigma)$-rule} $e = (\alpha, \gamma)$ associates,
to any model $R$ of $(S,\Sigma)$,
a pair $(\alpha^R, \gamma^R)$ of \emph{parallel} 
morphisms in a suitable category.
In a sense made precise later, we can ask whether 
  \[  \alpha^R \enspace \leq \enspace \gamma^R \enspace , \]
due to our use of relative monads towards families of \emph{preordered} sets.
If this is the case, we say that $R$ \emph{satisfies} the rule $e$.
A \emph{2--signature} is a pair $((S,\Sigma),A)$ consisting of a 1--signature $(S,\Sigma)$, which specifies the types and 
terms of a language, together with a family $A$ of 
$(S,\Sigma)$-rules, which specifies reduction rules on those terms.
Given a 2--signature $((S,\Sigma),A)$, we call a \emph{model of $((S,\Sigma),A)$} any model of $(S,\Sigma)$ that
satisfies each rule of $A$.
The \emph{category of models of $((S,\Sigma),A)$} is defined to be the full subcategory of models of $(S,\Sigma)$ 
whose objects are models of $((S,\Sigma),A)$.

We would like to exhibit an initial object in the category of models of $((S,\Sigma),A)$,
and thus must filter out rules which are never satisfied.
We call \emph{elementary} $(S,\Sigma)$-rule any $(S,\Sigma)$-rule whose codomain is of a particular form.
Our main result states that for any family $A$ of elementary $(S,\Sigma)$-rules, the category of models of $((S,\Sigma),A)$ has an initial object.
The class of elementary rules is large enough to account for the 
fundamental reduction rules; 
in particular, $\beta$- and $\eta$-reductions are given by elementary rules.
However, our notion of rule does only capture \emph{axioms}, that is, unconditional rewrite rules.

Our definitions ensure that any reduction rule between terms that is expressed by a directed equation $e\in A$ is automatically 
propagated into subterms in the initial model. 
The family $A$ of rules hence only needs to contain some \enquote{generating} rules, 
a fact that is well illustrated by the example 2--signature $\Lambda\beta$ of the untyped lambda calculus with $\beta$-reduction \cite{ahrens_relmonads}:
this signature has only one directed equation $\beta$ which expresses $\beta$-reduction at the root of a term,
   \[ (\lambda x . M)N \enspace \leadsto \enspace M [x:=N] \enspace .\]  
The initial model of $\Lambda\beta$ is given by the untyped lambda calculus, equipped with
the reflexive and transitive $\beta$-reduction relation $\twoheadrightarrow_{\beta}$ as presented 
in \cite{barendregt_barendsen}.

\subsection{Related Work}\label{sec:rel_work}

Initial semantics results for syntax with variable binding were first presented 
on the LICS'99 conference (see below).
Those results are concerned only with the \emph{syntactic aspect} of languages: they characterize the \emph{terms}
of a language as an initial object in some category, while not taking into account reductions on terms.
In lack of a better name, we refer to this kind of initiality results as \emph{purely syntactic}.

Some of these initiality theorems have been extended to also incorporate semantic aspects, e.g., in form of 
equivalence relations between terms. These extensions are reviewed further below.

\subsubsection{Purely syntactic results} \label{sec:purely_syntactic}

Initial Semantics for ``pure'' syntax---i.e.\ without considering semantic aspects---%
with variable binding were presented by several people independently, differing in the representation of 
variable binding:

The \emph{nominal approach} by Gabbay and Pitts \cite{gabbay_pitts99} (see also \cite{gabbay_pitts, pitts}) 
uses a set theory enriched with \emph{atoms} to establish an initiality result. 
Their approach models lambda abstraction as a constructor which takes a pair of a variable name and a term as arguments.
In contrast to the other techniques mentioned in this list, in the nominal approach
syntactic equality is different from $\alpha$--equivalence.
Hofmann \cite{hofmann} proves an initiality result modelling variable binding in a 
Higher--Order Abstract Syntax (HOAS) style.
Fiore, Plotkin, and Turi \cite{fpt} (see also subsequent work by Fiore \cite{fio02, fio05}) 
model variable binding through nested datatypes as introduced in \cite{BirdMeertens98:Nested}.
The approach of \cite{fpt} is extended to \emph{simply--typed syntax} in \cite{DBLP:conf/ppdp/MiculanS03}.
In \cite{Tanaka:2005:UCF:1088454.1088457} the authors
generalize and subsume those three approaches to a general category of contexts. 
An overview of this work and references to more technical papers is given in \cite{Power:2007:ASS:1230146.1230276}.
In \cite{DBLP:conf/wollic/HirschowitzM07} the authors prove an initiality result for untyped syntax based on 
the notion of \emph{module over a monad}. Their work is extended to simply--typed syntax in \cite{ju_phd}.

\subsubsection{Incorporating Semantics}

Rewriting in nominal settings has been examined in \cite{fernandez_gabbay_nominal_rewriting}.
In \cite{DBLP:journals/njc/GhaniL03} the authors present rewriting for algebraic theories without variable binding;
they characterize 
equational theories (with a \emph{symmetry} rule) resp.\ rewrite systems (with \emph{reflexivity} and \emph{transitivity} rule, but without \emph{symmetry})
 as \emph{coequalizers} resp.\ \emph{coinserters} in a category of monads on
the categories $\Set$ of sets resp.\ $\PO$ of preordered sets. 
Fiore and Hur \cite{DBLP:conf/icalp/FioreH07} have extended Fiore's work to integrate semantic aspects into initiality results.
In particular, Hur's thesis \cite{hur_phd} is dedicated to \emph{equational} systems for syntax with variable binding.
In a ``Further research'' section \cite[Chap.\ 9.3]{hur_phd}, Hur suggests the use of relations 
to model ``rewrite systems'' with \emph{directed} equations.
Hirschowitz and Maggesi \cite{DBLP:conf/wollic/HirschowitzM07} prove initiality of the set of lambda terms modulo $\beta$- and $\eta$-conversion
in a category of \emph{exponential monads}.
In an unpublished paper \cite{journals/corr/abs-0704-2900}, the authors  define a notion of \emph{half--equation} and 
\emph{equation} to express congruence between terms. We adopt their definition in this paper, but interpret a pair
of half--equations as \emph{directed} equation rather than equation.
In a ``Future Work'' section \cite[Sect.\ 8]{DBLP:journals/iandc/HirschowitzM10}, they 
mention the idea of using preorders as an approach to model  semantics, 
and they suggest interpreting the untyped lambda calculus with $\beta$- and $\eta$-reduction rule as a monad over the category $\PO$ of preordered sets.
The present work gives an alternative viewpoint to their suggestion by considering
the lambda calculus with $\beta$-reduction---and a class of programming languages in general---%
as a preorder--valued \emph{relative} monad on the functor $\Delta:\Set\to\PO$.
The rationale underlying our use of relative monads from sets to preorders is that we consider \emph{contexts} 
to be given by unstructured sets, whereas \emph{terms} of a language carry
structure in form of a reduction relation. In this view it is reasonable to suppose variables and terms to live in \emph{different} 
categories, which is possible through the use of relative monads on the functor $\Delta:\Set\to\PO$ (cf.\ Definition\ \ref{def:delta}) 
instead of traditional monads (cf.\ also \cite{ahrens_relmonads}).
Relative monads were introduced in \cite{DBLP:journals/corr/AltenkirchCU14}. In that work, the authors
characterize the untyped lambda calculus as a relative monad over the inclusion functor from finite sets to sets (see Example\ \ref{ex:ulc_finite_contexts}).
T.\ Hirschowitz \cite{DBLP:journals/corr/Hirschowitz13}
defines a category \fixedtextsf{Sig} 
of 2--signatures for \emph{simply--typed} syntax with reduction rules, and constructs an adjunction
between \fixedtextsf{Sig} and the category $\mathsf{2CCCat}$
of small cartesian closed 2--categories. He thus associates to any signature a 2--category of 
types, terms and reductions satisfying a universal property. 
More precisely, terms are given by morphisms in this category, and 
reductions are expressed by 2--cells between terms.
His approach differs from ours in the way in which variable binding is modelled:
Hirschowitz encodes binding in a Higher--Order Abstract Syntax (HOAS) style through exponentials.

\subsection{Synopsis}
In Section~\ref{sec:monads_modules} we review the definition of relative monads and modules over such monads as well as their morphisms.  
Some constructions on monads and modules are given, which will be of importance in what follows.

In Section~\ref{sec:comp_types_sem} we define arities, half--equations and directed equations, as well as their models.
Afterwards we prove our main result. 

In Section~\ref{sec:trans_pcf_ulc} we describe the formalization in the proof assistant \coq of an instance of our 
main result, for the particular case of the language $\PCF$.

Some conclusions and future work are stated in Section~\ref{sec:conclusion}.

\section{Relative monads and modules thereon}\label{sec:monads_modules}

Monads capture the notion of \emph{simultaneous substitution} and its properties.
This was first observed in \cite{alt_reus}, who characterize the lambda calculus and substitution on it
as a monad on the category of sets.

The functor underlying a monad is necessarily endo---this is enforced by the type of monadic multiplication.
\emph{Relative monads} were introduced in \cite{DBLP:journals/corr/AltenkirchCU14} to overcome this
restriction. One of their motivations was to consider the untyped lambda calculus over \emph{finite} contexts 
and a suitable substitution operation on it as 
a monad--like structure---in the spirit of aforementioned earlier work \cite{alt_reus}, but now with different categories for source and target.

We review the definition of relative monads and define suitable \emph{colax morphisms of relative monads}.
Afterwards we define \emph{modules over relative monads} and generalize the constructions on 
modules over monads defined in \cite{DBLP:conf/wollic/HirschowitzM07} to modules over \emph{relative} monads.

\subsection{Definitions}\label{sec:rel_mon_defs}

We review the definition of relative monad as given in 
\cite{DBLP:journals/corr/AltenkirchCU14}
and define suitable morphisms for them.
As an example we consider the lambda calculus with $\beta$-reduction as a relative monad from sets to preorders, on the functor 
$\Delta : \Set\to\PO$ (cf.\ Definition\ \ref{def:delta}).
Afterwards we define \emph{modules over relative monads} and carry over the constructions on modules over regular monads
of \cite{DBLP:conf/wollic/HirschowitzM07} to modules over relative monads.

The definition of relative monads is analogous to that of monads \emph{in Kleisli form},
except that the underlying map of objects is between \emph{different} categories.
Thus, for the operations to remain well--typed, one needs an additional \enquote{mediating} functor, in the following 
usually called $F$, which is inserted wherever necessary:

\begin{defi}[Relative monad, \cite{DBLP:journals/corr/AltenkirchCU14}]\label{def:relative_monad}
Given categories $\C$ and $\D$ and a functor $F : \C\to \D$, a \emph{relative monad $P : \C \stackrel{F}{\to}\D$ on $F$}
is given by the following data:
\begin{itemize}
 \item a map $P\colon \C\to\D$ on the objects of $\C$,
 \item 
for each object $c$ of $\C$, a morphism $\eta_c\in \D(Fc,Pc)$ and
 \item 
for each two objects $c,d$ of $\C$, a \emph{substitution} map (whose subscripts we usually omit)
 \[	\sigma_{c,d}\colon \D (Fc,Pd) \to \D(Pc,Pd)\]
\end{itemize}
such that the following diagrams commute for all suitable morphisms $f$ and $g$:
\begin{equation*}
\begin{xy}
\xymatrix @=4pc{
Fc \ar [r] ^ {\we_c} \ar[rd]_{f} & Pc \ar[d]^{\kl{f}} & Pc \ar@/^1pc/[rd]^{\kl{\we_c}} \ar@/_1pc/[rd]_{\id}& {} & 
 Pc \ar[r]^{\kl{f}} \ar[rd]_{\kl{\comp{f}{\kl{g}}}} & Pd\ar[d]^{\kl{g}} \\
{} & Pd \enspace , & {} & Pc \enspace , & {} & Pe \enspace .\\
}
\end{xy}
\end{equation*}
We sometimes omit the adjective \enquote{relative} when it is clear from the context that the monad
referred to is a relative monad.
\end{defi}

\noindent
\begin{rem}
 Relative monads on the identity functor $\Id:\C\to\C$ precisely correspond to monads.
\end{rem}

Various examples of relative monads are given in \cite{DBLP:journals/corr/AltenkirchCU14}.
They give one example related to syntax and substitution:

\begin{exa}[Lambda calculus over finite contexts]\label{ex:ulc_finite_contexts}
 In \cite{DBLP:journals/corr/AltenkirchCU14}, the authors consider the untyped lambda calculus
 as a relative monad on the functor $J:\Fin_{\mathrm{skel}} \to\Set$.
 Here the category $\Fin_{\mathrm{skel}}$ is the category of finite cardinals, i.e.\ 
 the skeleton of the category $\Fin$ of \emph{finite}
  sets and maps between finite sets.
 The category $\Set$ is the category of sets, cf.\ Definition\ \ref{def:SET}.
\end{exa}

We will give another example (cf.\ Example\ \ref{ex:ulcbeta}) of how to view syntax \emph{with reduction rules} as a relative monad.
For this, we need some more definitions.

\begin{defi}\label{def:SET}
  The category $\Set$ is the category of sets and total maps between them, together with 
  the usual composition of maps.
\end{defi}
\begin{defi}\label{def:cat_PO}
 The category $\PO$ of preorders has, as objects, sets equipped with a preorder, and, as morphisms between any two 
 preordered sets $A$ and $B$, the monotone functions from $A$ to $B$.
The category $\PO$ is cartesian closed:
 given $f, g\in \PO(A,B)$, we say that $f \leq g$ if and only if
     for any $a\in A$, $f(a) \leq g(a)$ in $B$. 
\end{defi}

\begin{defi}[Functor $\Delta:\Set\to\PO$ and forgetful functor]
\label{def:delta}
\label{def:adj_set_po}
 We call $\Delta : \SET\to\PO$ the left adjoint of the forgetful functor $U : \PO\to\Set$,

\begin{equation*}   
       \begin{xy}
        \xymatrix @C=4pc {
                  **[l]\Set \rtwocell<5>_U^{\Delta}{'\bot} & **[r]\PO
}
       \end{xy} \enspace .
\end{equation*}

\noindent
The functor $\Delta$ associates,  
 to each set $X$, the set itself together with the smallest preorder, i.e.\ the diagonal of $X$,
\[ \Delta(X) := (X,\delta_X) \enspace . \] 
 In other words, for any $x,y\in X$ we have $x\delta_X y$ if and only if $x = y$.
 The functor $\Delta:\Set\to \PO$ is a \emph{full embedding}, i.e.\ it is fully faithful and injective on objects.
 We have $\comp{\Delta}{U} = \Id_{\Set}$. Altogether,
the embedding $\Delta:\SET\to\PO$ is a coreflection.
We denote by $\varphi$ the family of isomorphisms
   \[   \varphi_{X,Y} : \PO(\Delta X, Y) \cong \Set(X,UY) \enspace . \]
We omit the indices of $\varphi$ whenever they can be deduced from the context.

\end{defi}

\begin{defi}[Category of families]\label{def:TST}\label{def:TS}
 Let $\C$ be a category and $T$ be a set, i.e., a discrete category.
 We denote by $\family{\C}{T}$ the functor category, an object of which is a $T$--indexed family of objects of $\C$.  
 Given two families $V$ and $W$, a morphism $f : V \to W$ is a family of morphisms in $\C$,
  \[ f : t \mapsto f(t) : V(t) \to W(t) \enspace . \]
We write $V_t := V(t)$ for objects and morphisms.
 Given another category $\D$ and a functor $F : \C\to \D$, we denote by $\family{F}{T} : \family{\C}{T} \to \family{\D}{T}$ the functor
  given by postcomposition with $F$, that is, $f \mapsto \comp{f}{F}$.
  
\end{defi}
\begin{rem}\label{rem:adj_set_po_typed}
Note that an adjunction $L\dashv R : \C \to \D$ carries over to functor categories
  $L^{\X} \dashv R^{\X} : \C^\X \to \D^\X$.
 
We use this fact in the following case: given a set $T$, the adjunction of Definition\ \ref{def:adj_set_po} induces an adjunction
\[
   \begin{xy}
        \xymatrix@C=4pc{
                   **[l]\family{\Set}{T}\rtwocell<5>_{\family{U}{T}}^{\family{\Delta}{T}}{'\bot} & **[r]\family{\PO}{T}
}
       \end{xy} \enspace .
\]

\end{rem}
 
\begin{exa}[Simply--typed lambda calculus as relative monad on $\TDelta{T}$]\label{ex:ulcbeta}
In this example, we consider $\TDelta{T} : \TS{T} \to \TP{T}$. In the \coq code shown here, \lstinline!Type! plays the role of $\Set$.
Let \[T := \TLCTYPE ::= ~ * ~ \mid ~\TLCTYPE \TLCar \TLCTYPE\] be the set of types of the simply--typed lambda calculus.
Consider the set family of simply--typed lambda terms over $\TLCTYPE$, indexed by typed contexts:
\begin{lstlisting}
Inductive TLC (V : T -> Type) : T -> Type :=
  | Var : forall t, V t -> TLC V t
  | Abs : forall s t TLC (V + s) t -> TLC V (s ~> t)
  | App : forall s t, TLC V (s ~> t) -> TLC V s -> TLC V t.
\end{lstlisting}
Here the context \lstinline!V + s! is the context \lstinline!V! extended by a fresh variable of type \lstinline!s!---the variable 
that is bound by the constructor \lstinline!Abs! (cf.\ also Section\ \ref{sec:deriv_and_fibre}).
We leave the object type arguments implicit and
write $\lambda M$ and $MN$ for $\Abs M$ and $\App M N$, respectively.
We equip each set $\TLC(V)(t)$ of lambda terms over context $V$ of object type $t$  with a preorder taken as the reflexive--transitive closure of the 
relation generated by the rule
\begin{equation*} \quad (\lambda M)N ~ \leq ~ M [*:= N] \end{equation*}
and its propagation into subterms.
 This defines a monad $\TLC_{\beta}$ from families of sets to families of preorders over the functor $\TDelta{T}$,
    \[\TLCB : \TS{T}\stackrel{\TDelta{T}}{\longrightarrow}\TP{T} .\] 
 The family $\we^{\TLC}$ is given by the constructor $\Var$, and the substitution map 
  \begin{equation}
     \sigma_{X,Y} : \TP{T}\bigl(\TDelta{T} (X), \TLCB(Y)\bigr) \to \TP{T}\bigl(\TLCB(X),\TLCB(Y)\bigr) \label{eq:subst}
  \end{equation}
 is given by capture--avoiding simultaneous substitution. 
Via the adjunction of Remark\ \ref{rem:adj_set_po_typed} the substitution can also be read as
  \[\sigma_{X,Y} : \TS{T}\bigl(X, \TLC(Y)\bigr) \to \TP{T}\bigl(\TLCB(X),\TLCB(Y)\bigr) \enspace .\]

\end{exa}

\noindent
In the previous example, the substitution of the lambda calculus satisfies an additional monotonicity property:
the map $\sigma_{X,Y}$ in Display\ \eqref{eq:subst} is monotone for the preorders on hom--sets defined in Definition\ \ref{def:cat_PO}.
This motivates the following definition:

\begin{defi}\label{def:enriched_monad}
  Let $P$ be a relative monad on $\TDelta{T}$ for some set $T$. We say that $P$ is \emph{a reduction monad} if, for any $X$ and $Y$,
  the substitution $\sigma_{X,Y}$ is \emph{monotone} for the preorders on $\family{\PO}{T}(\TDelta{T} X, PY)$
   and $\family{\PO}{T}(PX, PY)$.
\end{defi}
The monad $\TLCB$ of Example\ \ref{ex:ulcbeta} is thus a reduction monad.

\begin{rem}[Relative monads are functorial]
  \label{rem:rel_monad_functorial}
Given a monad $P$ over $F : \C\to\D$, a functorial action
for $P$ is defined by setting, for any morphism $f:c\to d$ in $\C$,
\[
   P(f) := \lift_P(f) := \kl{\comp{Ff}{\we}} \enspace.
\]
The functor axioms are easily proved from the monadic axioms.
\end{rem}

The substitution $\sigma = (\sigma_{c,d})_{c,d \in |\C|}$ of a relative monad $P$  is binatural: 
\begin{rem}[Naturality of substitution]\label{rem:subst_natural_rel}
Given a relative monad $P$ over $F:\C \to \D$, then its substitution $\sigma$ is natural in $c$ and $d$.
We write $f^*(h) := \comp{f}{h}$.
 For naturality in $c$ we check that the diagram 
 \[
  \begin{xy}
   \xymatrix{
     c \ar[d]_{f}  &    **[l]\D(Fc,Pd) \ar[r]^{\sigma_{c,d}}  &  **[r]\D (Pc,Pd)\\
     c'             &    **[l]\D(Fc',Pd) \ar[r]_{\sigma_{c',d}} \ar[u]^{(Ff)^*}  &  **[r]\D (Pc',Pd) \ar[u]_{(Pf)^*}
}
  \end{xy}
 \]
commutes. 
Given $g\in \D(Fc',Pd)$, we have 
\begin{align*} 
        \comp{Pf}{\sigma(g)} &= \comp{\sigma(\comp{Ff}{\eta_{c'}})}{\sigma(g)} \\
                             &\stackrel{3}{=} \sigma(\comp{\comp{Ff}{\eta_{c'}}}{\sigma(g)}) \\ 
                             &\stackrel{1}{=} \sigma(\comp{Ff}{g})  \enspace ,   
\end{align*}
where the numbers correspond to the diagrams of Definition\ \ref{def:relative_monad} used
to rewrite in the respective step.
Similarly we check naturality in $d$. Writing $h_*(g) := \comp{g}{h}$, the diagram

 \[
  \begin{xy}
   \xymatrix{
     d \ar[d]_{h}  &    **[l]\D(Fc,Pd) \ar[r]^{\sigma_{c,d}} \ar[d]_{(Ph)_*} &  **[r]\D (Pc,Pd)\ar[d]^{(Ph)_*} \\
     d'             &    **[l]\D(Fc',Pd) \ar[r]_{\sigma_{c',d}}   &  **[r]\D (Pc',Pd) 
}
  \end{xy}
 \]
commutes: given $g\in \D(Fc, Pd)$, we have
\begin{align*} 
        \comp{\sigma(g)}{Ph} &= \comp{\sigma(g)}{\sigma(\comp{Fh}{\eta_{d'}})} \\
                             &\stackrel{3}{=} \sigma (\comp{g}{\sigma(\comp{Fh}{\eta_{d'}})})\\
                             &= \sigma(\comp{g}{Ph})  \enspace .   
\end{align*}
\end{rem}

\noindent
If $(T_i)_{i\in I}$ is a family of sets and $f : I \to J$ a map of sets, then we obtain a 
family of sets $(T'_j)_{j \in J}$ by setting $T'_j := \coprod_{\{i | f(i) = j\} } T_i$.
The following construction generalizes this reparametrization:

\begin{defi}[Retyping functor]\label{def:retyping_adjunction_kan}
Let $T$ and $T'$ be sets and $g:T\to T'$ be a map.
Let $\C$ be a cocomplete category.
The map $g$ induces a functor 
 \[ g^*:\family{\C}{T'} \to \family{\C}{T} \enspace, \quad W \mapsto \comp{g}{W} \enspace . \]
  The \emph{retyping functor associated to $g:T\to T'$},
 \[ \retyping{g}:\family{\C}{T} \to \family{\C}{T'} \enspace,  \]
  is defined as the left Kan extension operation 
  along $g$, that is, we have an adjunction
\begin{equation}
  \begin{xy}
   \xymatrix @C=4pc {
            **[l]\family{\C}{T} \rtwocell<5>_{g^*}^{\retyping{g}}{'\bot} &  **[r]\family{\C}{T'}
}  
  \end{xy} \enspace .
\label{eq:adjunction_retyping}
 \end{equation}
  Put differently, the map $g : T\to T'$ induces an endofunctor $\bar{g}$ on $\family{\C}{T}$ with object map
 \[ \bar{g}(V) := \comp{g}{\retyping{g}(V)} \]
and we have a natural transformation \lstinline!ctype!, the unit of the adjunction of Display\ \eqref{eq:adjunction_retyping},
 \[  \text{\lstinline!ctype!}:\Id \nt \bar{g} : \family{\C}{T} \to \family{\C}{T} \enspace . \]
\end{defi}

\begin{defi}[Pointed index sets]\label{def:cat_indexed_pointed}\label{def:cat_set_pointed}
  Given a category $\C$, a set $T$ and a natural number $n\in \NN$, we denote by $\family{\C}{T}_n$ the category
  with, as objects, diagrams of the form
   \[ n \stackrel{\vectorletter{t}}{\to} T \stackrel{V}{\to} \C \enspace , \]
  written $(V, t_1, \ldots, t_n)$ with $t_i := \vectorletter{t}(i)$.
  A morphism $h$ to another such $(W,\vectorletter{t}) $
  with the same pointing map $\vectorletter{t}$ is given by a morphism $h : V\to W$ in $\family{\C}{T}$;
  there are no morphisms between objects with different pointing maps.
  Any functor $F : \family{\C}{T} \to \family{\D}{T}$ extends to $F_n : \family{\C}{T}_n \to \family{\D}{T}_n$ via
   \[ F_n (V,t_1,\ldots,t_n) := (FV, t_1,\ldots,t_n) \enspace . \]
\end{defi}

\begin{rem}
 The category $\family{\C}{T}_n$ consists of $T^n$ copies of $\family{\C}{T}$, which do not interact.
 Due to the ``markers'' $(t_1, \ldots, t_n)$ we can act differently on each copy, 
  cf., e.g., Definitions\ \ref{def:derived_rel_mod_II} and \ref{def:fibre_rel_mod_II}.
\end{rem}

Retyping functors generalize to categories with pointed indexing sets;
when changing types according to a map of types $g:T\to T'$, the markers must be adapted as well:

\begin{defi}\label{def:retyping_functor_pointed}
Given a map of sets $g:T\to T'$, by postcomposing the pointing map with $g$, the retyping functor generalizes to the functor
 \[ \retyping{g}(n) : \family{\C}{T}_n \to \family{\C}{T'}_n \enspace , \quad (V, \vectorletter{t}) \mapsto \bigl(\retyping{g} V, g_*(\vectorletter{t})\bigr) \enspace ,  \] 
 where $g_*(\vectorletter{t}) := \comp{\vectorletter{t}}{g} : n\to T'$.
\end{defi}

Finally there is also a category where families of objects of $\C$ over different indexing sets are mixed together:

\begin{defi}\label{def:cat_TEns}
 Given a category $\C$, we denote by $\T\C$ the category where an object is a pair $(T,V)$ of 
a set $T$ and a family $V\in \family{\C}{T}$ of objects of $\C$ indexed by $T$.
   A morphism $(g,h)$ to another such $(T',W)$ is given by a map $g : T\to T'$ and a 
  morphism $h : V\to\comp{g}{W}$ in $\family{\C}{T}$, that is, a
  family of morphisms in $\C$, indexed by $T$,
  \[ h_t : V_t \to W_{g(t)} \enspace . \]
Suppose $\C$ has an initial object, denoted by $0_{\C}$. 
Given $n\in \mathbb{N}$, we call $\hat{n} = (n, k\mapsto 0_\C)$ the object of $\T\C$ that associates to any 
$1\leq k \leq n$ the initial object of $\C$.
 We call $\T\C_n$ the slice category $\hat{n} \downarrow \T\C$.
An object of this category consists of an object $(T,V) \in\T\C$ whose indexing set ``of types'' $T$ is 
pointed $n$ times, written $(T,V,\vectorletter{t})$, where $\vectorletter{t}$ is a vector of
elements of $T$ of length $n$. A morphism $(g,h) :(T,V,\vectorletter{t}) \to (T',V',\vectorletter{t'})$
is a morphism $(g,h) :(T,V) \to (T',V')$ as above, such that $\vectorletter{t'} = \comp{g}{\vectorletter{t}}$.

We call 
$\T U_n : \T\C_n \to \Set$
the forgetful functor associating to any pointed family $(T,V, t_1,\ldots,t_n)$ the indexing set $T$.
Note that for a fixed set $T$, the category $\family{\C}{T}_n$ (cf.\ Definition\ \ref{def:cat_indexed_pointed}) 
is the fibre over $T$ of this functor.
\end{defi}

\begin{rem}[Picking out sorts]\label{rem:nat_trans_picking_sort}
   Let $1: \T\C_n \to \Set$ denote the constant functor which maps objects to the terminal object
   of the category $\Set$.
   A natural transformation $\tau:1 \nt \T U_n$ associates to any object $(T,V,\mathbf{t})$ of the 
   category $\T\C_n$
   an element of $T$. 
   Naturality imposes that $\tau(T',V',\mathbf{t'}) = g \left(\tau(T,V,\mathbf{t})\right)$ for any 
     $(g,h) : (T,V,\mathbf{t}) \to (T',V',\mathbf{t'})$.
\end{rem}
\begin{nota}\label{not:tau_simpl_notation}
 Given a natural transformation $\tau : 1 \nt \T U_n$ as in Remark\ \ref{rem:nat_trans_picking_sort}, we write
  \[ \tau (T,V,\mathbf{t}) := \tau (T,V,\mathbf{t})(*) \in T \enspace , \]
  i.e., we omit the argument $*\in 1_{\Set}$ of the singleton set.
\end{nota}

\begin{exa}
 For $ 1 \leq k \leq n$, we denote by $k : 1 \nt \T U_n : \T\C_n \to \Set$ the natural transformation such that $k(T,V,\vectorletter{t}) := \vectorletter{t}(k)$.
\end{exa}

We are interested in monads on the category $\TS{T}$ of families of sets indexed by $T$ and relative monads on 
$\TDelta{T}:\TS{T}\to\TP{T}$ as well as their relationship:

\begin{defi}[Relative monads on $\TDelta{T}$ and monads on $\TS{T}$]\label{def:rmon_delta_endomon}
  Let $P$ be a relative monad on $\TDelta{T}:\TS{T}\to\TP{T}$. 
 By postcomposing with the forgetful functor $\family{U}{T}:\TP{T}\to \TS{T}$ we obtain a monad 
   \[UP : \TS{T}\to \TS{T}\enspace . \]
   The substitution is defined, for $m : X \to UPY$ by setting
   \[U\sigma : m\mapsto  U \left(\kl{\varphi^{-1} m}\right) \enspace ,  \]
making use of the adjunction $\varphi$ of Remark\ \ref{rem:adj_set_po_typed}.
Conversely, to any monad $P$ over $\TS{T}$, we associate a relative monad $\Delta P$ over $\TDelta{T}$ by 
 postcomposing with $\TDelta{T}$.
\end{defi}

\noindent
We generalize the definition of colax monad morphisms \cite{Leinster2003} to relative monads:

\begin{defi}[Colax morphism of relative monads] 
 \label{def:colax_rel_mon_mor}
Suppose given two relative monads $P : \C\stackrel{F}{\to}\D$ and 
$Q : \C'\stackrel{F'}{\to}\D'$. A \emph{colax morphism of relative monads} from $P$ to $Q$ is a
quadruple $h = (G,G',N,\tau)$ of
a functor $G\colon \C\to\C'$, a functor $G': \D \to \D'$ as well as a natural transformation
$N: F'G \to G'F$ and a natural transformation $\tau : PG' \nt GQ$
such that the following diagrams commute for any objects $c,d$ and any suitable morphism $f$:
\begin{equation*} 
 \begin{xy}
  \xymatrix @=5pc{
  G'Pc \ar[r]^{G'\kl[P]{f}} \ar[d]_{\tau_c}& G'Pd \ar[d]^{\tau_d} \\
  QGc \ar[r]_{\kl[Q]{\comp{Nc}{\comp{G'f}{\tau_d}}}} & QGd 
}
 \end{xy}
\qquad
\begin{xy}
  \xymatrix @=5pc {
 F'Gc \ar[r]^{Nc} \ar[rrd]_{\we^Q_{Gc}} & 
                 G'F c \ar[r]^{G'\we^P_c}& G'Pc \ar[d]^{\tau_c} \\
{} & {} & QGc.
}
\end{xy}
\end{equation*}
 Naturality of the family $\tau$ is actually a consequence of the commutative diagrams and 
 may be omitted from the definition.
\end{defi}

\begin{rem} \label{rem:rel_mon_mor_case}
  In Section\ \ref{sec:comp_types_sem} we are going to use the following instance of the preceding definition:
  the categories $\C$ and $\C'$ are instantiated by $\TS{T}$ and $\TS{T'}$, respectively, for sets $T$ and $T'$.
  The functor $G$ is the retyping functor (cf.\ Definition\ \ref{def:retyping_adjunction_kan}) associated to some 
  translation of types $g : T \to T'$.
  Similarly, the categories $\D$ and $\D'$ are instantiated by $\TP{T}$ and $\TP{T'}$, and 
  the functor $F$ by
  \[F := \TDelta{T} : \TS{T} \to \TP{T} \enspace ,\]
  and similar for $F'$:
\[
 \begin{xy}
  \xymatrix@!=2.5pc{
       **[l] \TS{T} \ar[r]^{\TDelta{T}} \ar[d]_{\retyping{g}}   & **[r] \TP{T} \ar[d]^{\retyping{g}}\\
       **[l] \TS{T'}\ar[r]_{\TDelta{T'}} & **[r]\TP{T'}. \ultwocell<\omit>{\Id}
}
 \end{xy}
\]
\end{rem}

Given a monad $P$ on $F:\C\to\D$, the notion of \emph{module over $P$} generalizes the notion of monadic substitution:

\begin{defi}[Module over a relative monad]
 \label{def:rmodule}
Let $P\colon\C\stackrel{F}{\to}\D$ be a relative monad and let $\E$ be a category. A \emph{module $M$ over $P$ with codomain $\E$} is given by
\begin{itemize}
 \item a map $M: \C \to \E$ on the objects of the categories involved and 
 \item for all objects $c,d$ of $\C$, a map 
      \[ 
          \varsigma_{c,d} : \D (Fc,Pd) \to \E (Mc,Md)
      \]
\end{itemize}
such that the following diagrams commute for all suitable morphisms $f$ and $g$:
\begin{equation*}
\begin{xy}
\xymatrix @=3pc{
 Mc \ar[r]^{\mkl{f}} \ar[rd]_{\mkl{\comp{f}{\kl{g}}}} & Md\ar[d]^{\mkl{g}} & Mc  \ar@/^1pc/[rd]^{\mkl{\we_c}} \ar@/_1pc/[rd]_{\id} & {} \\
  {} & Me & {} & Mc. \\
}
\end{xy}
\end{equation*}
\end{defi}

\noindent
A functorial action for such a module $M$ is then defined similarly to that for relative monads:
 for any morphism $f : c \to d$ in $\C$ we set
 \[ M(f) := \rmlift_M(f) := \varsigma({ \comp{Ff}{\eta}}) \enspace .\]

Examples of modules over relative monads are given in the next section.

A \emph{module morphism} is a family of morphisms that is compatible with module substitution in the source and target modules:
\begin{defi}[Morphism of relative modules]\label{def:rel_mod_mor}
 Let $M$ and $N$ be two relative modules over $P\colon\C\stackrel{F}{\to}\D$ with codomain $\E$. 
A \emph{morphism of relative $P$--modules} from $M$ to $N$ is given by a collection of morphisms $\rho_c\in\E(Mc,Nc)$ 
such that for all morphisms $f\in \D(Fc,Pd)$ the following diagram commutes:
\begin{equation*}
 \begin{xy}
  \xymatrix @=3pc{
  Mc \ar[r]^{\mkl[M]{f}} \ar[d]_{\rho_c}& Md \ar[d]^{\rho_d}  \\
  Nc \ar[r]_{\mkl[N]{f}} & Nd.\\
}
 \end{xy}
\end{equation*}
\end{defi}

\noindent
The modules over $P$ with codomain $\E$ and morphisms between them form a category called $\RMod{P}{\E}$. 
Composition and identity morphisms of modules are defined by pointwise composition and identity.

\subsection{Constructions on modules over relative monads}\label{subsection:rmod_examples}
The following constructions generalize those used in \cite{DBLP:conf/wollic/HirschowitzM07} to
modules over \emph{relative} monads.

Any relative monad $P$ comes with the \emph{tautological} module over $P$ itself:
\begin{defi} [Tautological module] 
Every relative monad $P$ on $F:\C\to\D$ yields a module $(P,\sigma^P)$ --- also denoted by $P$ --- over itself, 
i.e.\ an object in the category $\RMod{P}{\D}$. 
\end{defi}

\begin{exa}[Example\ \ref{ex:ulcbeta} continued]\label{ex:ulcb_taut_mod}
 The map $\TLCB : V \mapsto \TLCB(V)$ yields a module over the relative monad $\TLCB$, the \emph{tautological $\TLCB$--module} $\TLCB$.
\end{exa}

Constant functors are modules over relative monads with the same source category:

\begin{defi}[Constant and terminal module] 
Let $P$ be a relative monad on $F:\C\to\D$.
For any object $e \in \E$ the constant map $T_e\colon\C\to\E$, $c\mapsto e$ for all $c\in \C$, is equipped with the structure of a $P$--module
by setting $\varsigma_{c,d}(f) = \id_e$. 
In particular, if $\E$ has a terminal object $1_\E$, then the constant module $T_{1_\E} : c \mapsto 1_\E$ is terminal in $\RMod{P}{\E}$.
\end{defi}

Constant modules can be expressed via the following construction:

\begin{defi}[Postcomposition with a functor]
 Let $P$ be a relative monad on $F:\C\to\D$, and let $M$ be a $P$--module with codomain $\E$.
 Let $G: \E \to \X$ be a functor. Then the object map $\comp{M}{G}:\C\to \X$ defined by $c\mapsto G(M(c))$ 
 is equipped with a $P$--module structure by setting, for $c, d \in \C$ and $f\in \D(Fc,Pd)$,
  \[  \varsigma^{\comp{M}{G}}(f) := G(\varsigma^M(f)) \enspace . \]
 For $M:=P$ (considered as tautological module over itself) and $G$ a constant functor mapping to an object $x\in \X$ and its identity morphism $\id_x$, 
 we obtain the constant module $(T_x,\id)$ as in the preceding definition.
\end{defi}

Postcomposition with a functor $F : \E\to\X$ extends to morphisms of modules and yields a functor between categories of modules, 
$\RMod{P}{\E} \to \RMod{P}{\X}$. We omit the details since we do not use this fact.

The next construction on modules does not change the target category, but the underlying relative monad: given a module $N$ over a relative monad $Q$ and a monad morphism $\tau : P \to Q$ into $Q$, we rebase or ``pull back''
the module $N$ along $\tau$, yielding a module over $P$:

\begin{defi}[Pullback module] \label{def_pullback}
 Suppose given two relative monads $P$ and $Q$ and a morphism $h = (G, G', N, \tau) : P \to Q$ as in Definition\ \ref{def:colax_rel_mon_mor}.
Let $M$ be a $Q$-module with codomain $\E$. We define a $P$-module $h^* M$ to $\E$ with object map 
 \[c\mapsto M (Gc) \]
by defining the substitution map, for $f : Fc \to Pd$, as
 \[\mkl[h^*M] f := \mkl[M]{\comp{N_c}{\comp {G'f}{\tau_d}}} \enspace . \] 
The module thus defined is called the \emph{pullback module of $M$ along $h$}. 
The pullback extends to module morphisms and is functorial. 
\end{defi}

\begin{defi}[Module morphism induced by a monad morphism]
\label{def:ind_rmod_mor}
With the same notation as in Definition\ \ref{def_pullback}, the monad morphism $h$ induces a morphism of $P$--modules $h:G'P \to h^*Q$. %
Note that the domain module is the module obtained by postcomposing (the tautological module of) $P$ with $G'$, 
whereas for (traditional) monads the domain module was just the tautological module of the domain monad \cite{DBLP:conf/wollic/HirschowitzM07}.
\end{defi}

We can take the \emph{product} of two modules:

\begin{defi}[Products lift] 
Suppose the category $\E$ has products. Let $M$ and $N$ be $P$--modules with codomain $\E$. Then the map 
\[ M \times N : \C\to\E, \quad c \mapsto Mc \times Nc \] 
 is canonically equipped with a substitution and thus constitutes a module called the \emph{product of $M$ and $N$}. 
 This construction extends to a product on $\RMod{P}{\E}$. 
\end{defi}

\begin{exa}\label{ex:ulcb_prod_mod}
 Given $s,t\in T$, the map $V\mapsto \TLCB(V)(s \TLCar t)\times \TLCB(V)(s)$ inherits a structure of an $\TLCB$--module.
\end{exa}

\subsection{Derivation \& fibre} \label{sec:deriv_and_fibre}

We are particularly interested in relative monads on the functor $\TDelta{T}: \TS{T}\to\TP{T}$ for some set $T$, 
 and modules over such monads.
\emph{Derivation} and \emph{fibre}, two important constructions of \cite{DBLP:journals/iandc/HirschowitzM10} on modules over monads on families of sets,
carry over to  modules over relative monads on $\family{\Delta}{T}$.

\begin{defi}\label{def_ctxt_ext}
Given $u\in T$, we denote by $D(u)\in \TS{T}$ the context with $D(u)(u) = \{*\}$ and $D(u)(t) = \emptyset$ for $u \neq t$.
For a context $V \in \TS{T}$ we set $V^{u} := V + D(u)$.
\end{defi}

\begin{defi}
  \label{def:rel_module_deriv}
Given a monad $P$ over $\TDelta{T}$ and a $P$--module $M$ with codomain $\E$, 
we define the derived module of $M$ with respect to $u\in T$ by setting
\[ M^u(V) := M (V^{u}) \enspace . \]
The module substitution is defined, for $f\in \TP{T}(\TDelta{T}V, PW)$, by 
 \[ \mkl[M^u]{f} := \mkl[M]{\shift{f}{u}} \enspace .\]
Here the ``shifted'' map 
   \[ \shift{f}{u} \in \TP{T}\bigl(\TDelta{T}(V^{u}), P(W^{u})\bigr) \]
is the adjunct under the adjunction of Remark\ \ref{rem:adj_set_po_typed}  
of the coproduct map
 \[ \varphi (\shift{f}{u}) := [\comp{f}{P(\inl)}, \we(\inr(*))] : V^{*u} \to UP(W^{u}) \enspace , \]
where $[\inl, \inr] = \id : W^{u}\to W^{u}$.
Derivation is an endofunctor on the category of $P$--modules with codomain $\E$. 
\end{defi}

\begin{exa}\label{ex:ulcb_der_mod}
 Given $V \in \TS{T}$ and $s\in T$, we denote by $V^s$ the context $V$ enriched by an additional variable of type $s$
as in Definition\ \ref{def_ctxt_ext}.
The map $\TLCB^s : V \mapsto \TLCB(V^s)$ inherits the structure of a $\TLCB$--module from the
  tautological module $\TLCB$ (cf.\ Example\ \ref{ex:ulcb_taut_mod}).
We call $\TLCB^s$ the \emph{derived module with respect to $s\in T$} of the module $\TLCB$; cf.\ also Section\ \ref{subsection:rmod_examples}.
\end{exa}

\begin{nota}\label{not:deriv_rmod_untyped}
 In case the set $T$ of types is $T = \{*\}$ the singleton set of types, i.e., when talking about untyped syntax,
we denote by $M'$ the derived module of $M$.
  Given a natural number $n$, we denote by $M^n$ the module obtained by deriving $n$ times the module $M$.
\end{nota}

We derive more generally with respect to a natural transformation $\tau : 1 \nt \T U_n$ as in Definition \ref{def:cat_TEns}:
\begin{defi}[Derived module]\label{def:derived_rel_mod_II}
  Let $\tau : 1 \nt \T U_n$ be a natural transformation.
 Let $T$ be a set and $P$ be a relative monad on $\TDelta{T}_n$.
 \noindent
Given any $P$--module $M$, we call \emph{derivation of $M$ with respect to $\tau$} the module with object map
$M^{\tau} (V):= M\left(V^{\tau(V)}\right)$.
\end{defi}

\begin{defi}\label{def:fibre_rel_mod_II}
 Let $P$ be a relative monad over $F$, and $M$ a $P$--module with codomain $\E^T$ for some category $\E$.
 The \emph{fibre module $\fibre{M}{t}$ of $M$ with respect to $t\in T$} has object map
  \[ c \mapsto M(c)(t) = M(c)_t \]
  and substitution map
  \[ \mkl[{\fibre{M}{t}}]{f} := \bigl(\mkl[M]{f}\bigr)_t \enspace .  \]
\end{defi}

\begin{exa}
 Given $t \in T$, the map $V\mapsto \TLCB(V)(t) : \TS{T} \to \PO$ inherits a structure of a $\TLCB$--module, the 
\emph{fibre module $\fibre{\TLCB}{t}$ with respect to $t\in T$}.
\end{exa}

\noindent
This definition generalizes to fibres with respect to a natural transformation as in Definition\ \ref{def:derived_rel_mod_II}.

\begin{exa}[Examples\ \ref{ex:ulcb_taut_mod}, \ref{ex:ulcb_der_mod}, \ref{ex:ulcb_prod_mod} continued]\label{ex:ulcb_constructor_mod_mor}
 Abstraction and application are morphisms of $\TLCB$--modules:
 \begin{align*} \Abs_{s,t} &: \fibre{\TLCB^{s}}{t} \to \fibre{\TLCB}{s\TLCar t} \enspace , \\
                \App_{s,t} &: \fibre{\TLCB}{s\TLCar t} \times \fibre{\TLCB}{s} \to \fibre{\TLCB}{t} \enspace .
 \end{align*}
\end{exa}

The pullback operation commutes with products, derivations and fibres :

\begin{lem} \label{lem:rel_pb_prod}Let $\C$ and $\D$ be categories and $\E$ be a category with products. Let $P\colon \C\to \D$ and $Q\colon \C\to D$ be monads
over $F:\C\to\D$ and $F':\C'\to \D'$, resp., and $\rho : P \to Q$ a monad morphism. 
Let $M$ and $N$ be $P$--modules with codomain $\E$. The pullback functor is cartesian:
 \[ \rho^* (M \times N) \cong \rho^*M \times \rho^*N \enspace .\]
\end{lem}
\begin{lem} \label{lem:rel_pb_comm}
Consider the setting as in the preceding lemma, with $F = \TDelta{T}$, and $t\in T$.
 Then we have
\[ \rho^* (M^t) \cong (\rho^*M)^t \enspace . \]
This readily generalizes to general derivation as defined in Definition~\ref{def:derived_rel_mod_II},
\[\rho^* (M^\tau) \cong (\rho^*M)^\tau \enspace . \]
\end{lem}

\begin{lem}\label{lem:rel_pb_fibre}
  Suppose $N$ is a $Q$--module with codomain $\family{\E}{T}$, and $t\in T$. Then
 \[ \rho^*\fibre{M}{t} \cong \fibre{\rho^*M}{t} \enspace . \]
\end{lem}

\begin{defi}[Substitution of \emph{one} variable]\label{def:hat_P_subst_typed}
  Let $T$ be a (nonempty) set and let $P$ be a \emph{reduction monad} (cf.\ Definition\ \ref{def:enriched_monad}) 
  over $\family{\Delta}{T}$.
  For any $s,t\in T$ and $X\in \TS{T}$ we define a binary substitution operation
      \begin{align*} 
   \subst_{s,t}(X) : P(X^{s})_t \times P(X)_s &\to P(X)_t, \\ 
    (y,z)&\mapsto y [*:= z] := \kl{[\we_X , \lambda x.z] }(y) \enspace .
\end{align*}
Here $y:P(X^s)_t$ is of type $t$ and lives in a context $X^s$, which is $X$ extended by an object variable of object type $s$.
This object variable is substituted with $z:P(X)_s$, which is of type $s$ itself.
  For any pair $(s,t)\in T^2$, we thus obtain a morphism of $P$--modules
  \[ \subst^{P}_{s,t} : \fibre{{P}^{s}}{t} \times \fibre{{P}}{s} \to \fibre{{P}}{t} \enspace . \]
\end{defi}
\noindent
Observe that this substitution operation is  monotone in both arguments: monotonicity in the first argument
is a consequence of the monadic axioms. Monotonicity in the second argument
is a consequence of $P$ being a reduction monads (Definition\ \ref{def:enriched_monad}).

\section{Signatures, models, initiality}\label{sec:comp_types_sem}

We combine the techniques of earlier work \cite{ahrens_ext_init, ahrens_relmonads} in order
to obtain an initiality result for simple type systems with reductions on the term level.
As an example, we specify, via the iteration principle coming from the universal property, a 
semantically faithful translation 
from $\PCF$ with its usual reduction relation to the untyped lambda calculus with $\beta$-reduction.

More precisely, in this section we define a notion of 
\emph{2--signature} and suitable
 \emph{models} for such 2--signatures, such that the types and terms generated by the 2--signature, 
              equipped with reduction rules according to the directed equations specified by the 2--signature, 
              form the
 \emph{initial model}.
Analogously to our work on untyped syntax \cite{ahrens_relmonads}, 
we define a notion of \emph{2--signature} with two levels:
a \emph{syntactic} level specifying types and terms of a language, and, on top of that, a \emph{semantic} level
specifying reduction rules on the terms.

\subsection{1--Signatures}

A \emph{1--signature} specifies types and terms over these types. We give two presentations of 
1--signatures, a \emph{syntactic} one (cf.\ Definition\ \ref{def:n_comp_classic_arity_syn}) 
and a \emph{semantic} one (cf.\ Definition\ \ref{def:typed_sig_symantically}).
The syntactic presentation is the same as in earlier work \cite{ahrens_ext_init}.
The semantic presentation in the present work is adapted from \cite{ahrens_ext_init} to our use of \emph{relative} monads or,
to be more precise, \emph{reduction monads}, as compared to traditional monads used in \cite{ahrens_ext_init}.

\subsubsection{Signatures for Types}

We present \emph{type signatures}, which later are used to specify the \emph{object types} of the languages
we consider.
Such signatures and their models were first considered by Birkhoff \cite{birkhoff1935}. 
Intuitively, a type signature specifies the respective arities---i.e., the number of arguments---of 
a set of operations on a set.

\begin{defi}[Type signature]\label{def:raw_sig}
A \emph{type signature $S$} is a family of natural numbers, i.e., a set $J_S$ (of operation symbols) and a map 
(carrying the same name as the signature) $S : J_S\to \mathbb{N}$.
For $j\in J_S$ and $n\in \mathbb{N}$, we also write $j:n$ instead of $j \mapsto n$.
An element of $J$ resp.\ its image under $S$ is called an \emph{arity} of $S$.
\end{defi}
Intuitively, the meaning of $j:n$ is that $j$ represents an operation taking $n$ arguments.

\begin{exa}[Type signature of $\TLCTYPE$, Example\ \ref{ex:ulcbeta}]\label{ex:type_sig_SLC}
  The type signature of the types of the simply--typed lambda calculus is given by
  \[ S_{\SLC} := \{* : 0 \enspace , \quad (\TLCar) : 2 \}\enspace .\]
\end{exa}

To any type signature we associate a category of \emph{models}.
We call \emph{model of $S$} any set $U$ equipped with operations according to the signature $S$. %
A \emph{morphism of models} is a map between the underlying sets that is compatible with the 
operations on either side in a suitable sense.
Models and their morphisms form a category.
We give the formal definitions:

\begin{defi}[Model of a type signature]\label{def:rep_alg_sig}
A model $R$ of a type signature $S$ is given by 
\begin{itemize}
 \item a set $X$ and
 \item for each $j\in J_S$, an operation $j^R:X^{S(j)} \to X$.
\end{itemize}
In the following, given a model $R$, we write $R$ also for its underlying set. %
\end{defi}

\begin{defi}[Morphisms of models]\label{def:mor_raw_rep}
 Given two models $T$ and $U$ of the type signature $S$, a \emph{morphism} from $T$ to $U$ is 
  a map $f : T\to U$ such that, for any arity $j:n$ of $S$, we have
  \[  \comp{j^T}{f} = \comp{(\underbrace{f\times\ldots\times f}_{n \text{ times}})}{j^U} \enspace . \]
\end{defi}

\begin{exa}\label{ex:type_PCF}
  The language $\PCF$ \cite{Plotkin1977223, Hyland00onfull} 
is a simply--typed lambda calculus with a fixed point operator
  and arithmetic constants.
Let $J:= \{\Nat, \Bool, (\PCFar)\}$. The signature of the types of \PCF~is given by the arities 
   \[S_{\PCF}:= \lbrace\Nat:0\enspace ,\quad \Bool: 0 \enspace,\quad (\PCFar): 2 \rbrace \enspace .\]
 A model $T$ of $S_{\PCF}$ is given by a set $T$ and three operations,
   \[
      \Nat^T : T\enspace, \quad \Bool^T : T \enspace , \quad  (\PCFar)^T : T\times T\to T \enspace .
   \]
 Given two models $T$ and $U$ of $S_{\PCF}$, a morphism from $T$ to $U$ is a map $f : T\to U$ between the underlying sets
 such that, for any $s,t\in T$,
\begin{align*}
      f(\Nat^T) &= \Nat^U \enspace ,  \\
      f(\Bool^T) &= \Bool^U \quad \text{ and}\\
     f(s \PCFar^T t) &= f(s) \PCFar^U f(t) \enspace .
\end{align*}
\end{exa}

\subsubsection{Signatures for terms}\label{sec:term_sigs}

Before giving our definition of signature for terms, we consider the example of the simply-typed
lambda calculus and describe our goals using this example.

Let $\TLCTYPE$ be the initial model of the signature for types of Example\ \ref{ex:type_sig_SLC}.
The signature for simply-typed lambda terms over those types may be given as follows:
\begin{equation}\{ \abs_{s,t} :=  \bigl[([s],t)\bigr] \to (s\TLCar t) \enspace , 
       \quad \app_{s,t} := \bigl[([],s\TLCar t),([],s)\bigr]\to t\}_{s,t\in\TLCTYPE} \enspace . 
 \label{eq:sig_tlc_simple}
\end{equation}
Intuitively, this means that we have two families of operations, $\abs$---abstraction---and $\app$---application.
The abstraction (specified by) $\abs_{s,t}$ takes one argument, specified by a list of length one,
where this argument is of type $t$ and lives in a context extended by a variable of type $s$.
The return type is $s \TLCar t$.

The parameters $s$ and $t$ range over the set $\TLCTYPE$ of types, the initial model of the 
signature for types from Example\ \ref{ex:type_sig_SLC}.
Our goal is to consider models of the simply--typed lambda calculus in monads over categories 
of the form $\TS{T}$ for \emph{any} set $T$---provided that $T$ is equipped with a model of the signature 
$S_{\SLC}$. It thus is more suitable to specify the signature of the simply--typed lambda calculus as follows:
\begin{equation}\{ \abs :=  \bigl[([1],2)\bigr] \to (1\TLCar 2) \enspace , 
       \quad \app := \bigl[([],1\TLCar 2),([],1)\bigr]\to 2\}\enspace . 
  \label{eq:sig_tlc_higher_order}
\end{equation}
For any model $T$ of $S_{\SLC}$, the variables $1$ and $2$
range over elements of $T$.
In this way the number of abstractions and applications depends on the model $T$ of $S_{\SLC}$:
intuitively, a 
model of the above signature of Display\ \eqref{eq:sig_tlc_higher_order} 
over a model  $T$ of $\TLCTYPE$
has $T^2$ 
abstractions and $T^2$ applications---one for each pair of elements of $T$.

To abstractly capture this kind of arities given in terms of natural numbers,
we introduce a notion of \emph{degree}---given by a natural number $n\in \NN$---, and term arities of degree $n$.
A result stated towards the end of this section, Lemma~\ref{lem:family_of_mods_cong_pointed_mod_relative},
explains how the ``grouping'' and ``ungrouping'' of such arities works.

We start by giving a very syntactic notion of arity for terms, in Definition~\ref{def:n_comp_classic_arity_syn}.
This notion is later related to a more semantic notion, given in Definition~\ref{def:term-arity-semantic}.

\begin{defi}[Type of degree $n$]
  For $n \geq 1$, we call \emph{types of $S$ of degree $n$} the elements of the set $S(n)$ 
 of types associated to the signature $S$ with free variables in the set $\{1,\ldots,n\}$.
  We set $S(0):= \hat{S}$.
  Formally, the set $S(n)$ may be obtained as the initial model of the type signature $S$ enriched by $n$ nullary arities.
\end{defi}

Types of degree $n$ are used to form elementary arities of degree $n$:

\begin{defi}[Elementary arity of degree $n$]\label{def:n_comp_classic_arity_syn}
 An elementary arity for terms over the signature $S$ for types of degree $n$ is of the form
\begin{equation} 
  \bigl[([t_{1,1},\ldots,t_{1,m_1}], t_1), \ldots, ([t_{k,1},\ldots,t_{k,m_k}], t_k)\bigr] \to t_0 \enspace , 
  \label{eq:syntactic_arity_higher_degree}
\end{equation}
  where $t_{i,j}, t_i \in S(n)$.
 More formally, an elementary  arity of degree $n$ over $S$ is a pair 
  consisting of an element $t_0\in S(n)$ and
  a list of pairs.  where each pair itself consists of a list $[t_{i,1},\ldots,t_{i,m_i}]$ of elements of $S(n)$ and 
  an element $t_i$ of $S(n)$.
\end{defi}

  An elementary arity of the form given in Display\ \eqref{eq:syntactic_arity_higher_degree} denotes a 
   constructor---or a family of constructors, for $n > 0$---whose output type is $t_0$,
  and whose $k$ inputs are terms of type $t_i$, respectively. In the input $i$, for $1\leq i \leq k$, variables of type according to the list
  $[t_{i,1},\ldots,t_{i,m_i}]$ are bound by the constructor.

\begin{exa}
 The prototypical examples of elementary arities of degree $2$ are abstraction and application of Display\ \eqref{eq:sig_tlc_higher_order}.
\end{exa}

Compared to our work on pure syntax \cite{ahrens_ext_init} 
we have to adapt the \emph{semantic} definition of signatures for terms, since we now work 
with reduction monads on $\TDelta{T}$ for some set $T$ instead of monads over families of sets.
\begin{defi}[Relative $S$--Monad] \label{def:s-rmon}
  Given a type signature $S$, the \emph{category $\SigRMon{S}$ of relative $S$--monads} 
 is defined as the category whose objects are pairs $(T,P)$ of
  a model $T$ of $S$ and a reduction monad 
\[P : \TS{T}\stackrel{\TDelta{T}}{\longrightarrow} \TP{T} \enspace .\]
  A morphism from $(T,P)$ to $(T', P')$ is a pair $(g, f)$ of a morphism of $S$--models $g : T\to T'$ and a 
    morphism of relative monads $f : P\to P'$ over the retyping functor $\retyping{g}$ as in Remark\ \ref{rem:rel_mon_mor_case}.

   Given $n\in \mathbb{N}$, we write $\SigRMon{S}_n$ for the category whose objects are pairs $(T,P)$ of a model $T$ of $S$ and 
  a reduction monad $P$ over $\TDelta{T}_n$. A morphism from $(T,P)$ to $(T', P')$ is a pair $(g, f)$ of a morphism of 
      $S$--model $g : T\to T'$ and a 
   monad morphism $f : P\to P'$ over the retyping functor $\retyping{g}_n$ defined in Definition\ \ref{def:retyping_functor_pointed}.

\end{defi}

In the following we need a category in which we gather modules over \emph{different} relative monads.
More precisely, an object is given by a pair of a relative monad and a module on it. 
A morphism between two such pairs is given by a pair of a monad morphism and a module morphism, where
for the latter we use the pullback operation on modules to obtain modules over the same relative monad:

\begin{defi}[Large category $\LRMod{n}{S}{\D}$ of modules]
  \label{def:lrmod_typed}
  Given a natural number $n\in \mathbb{N}$, a type signature $S$ and a category $\D$, 
  we call $\LRMod{n}{S}{\D}$ %
  the category an object of which is a pair $(P,M)$ of a relative $S$--monad $P \in \SigRMon{S}_n$ and a $P$--module with codomain $\D$.
  A morphism to another such $(Q,N)$ is a pair $(f, h)$ of a morphism of relative $S$--monads 
$f : P \to Q$ in $\SigRMon{S}_n$ and a morphism of relative modules $h : M \to f^*N$.
\end{defi}

We sometimes just write the module---i.e.\ the second---component of an object or morphism
 of the large category of modules.
Given $M\in \LRMod{n}{S}{\D}$, we thus write $M(V)$ or $M_V$ for the value of the module on the object $V$.

Let $S$ be a type signature.
A \emph{half-arity over $S$ of degree $n$} is a functor from relative $S$-monads to the category of large modules of degree $n$:

\begin{defi}[Half-Arity over $S$ (of degree $n$)]
  \label{def:half_arity_degree_semantic_typed}
 Given a type signature $S$ and $n\in \mathbb{N}$, we call \emph{half-arity over $S$ of degree $n$} a functor
  \[ \alpha : \SigRMon{S} \to \LRMod{n}{S}{\PO} \enspace . \]
 which is a section to the forgetful functor forgetting both the module and the \enquote{points}.
\end{defi}

Intuitively, a half-arity of degree $n$ associates, to any $S$-monad $P$, a module over $P_n$.

We restrict ourselves to a class of such functors, starting with the \emph{tautological} module:

\begin{defi}[Tautological module of degree $n$]\label{def:taut_mod_pointed}
  Given $n\in \mathbb{N}$, any relative monad $R$ over $\TDelta{T}$ induces a monad $R_n$ over $\TDelta{T}_n$ 
 with object map $(V, t_1,\ldots, t_n) \mapsto (RV, t_1,\ldots,t_n)$.
To any relative $S$--monad $R$ we associate
  the tautological module of $R_n$, 
  \[\Theta_n(R):= (R_n,R_n) \in \LRMod{n}{S}{\TP{T}_n} \enspace . \]
\end{defi}

\noindent
Furthermore, we use \emph{canonical natural transformations} (cf.\ Definition\ \ref{def:canonical_nat_trans}) to build \emph{elementary} half--arities; 
these transformations specify context extension (derivation) and
selection of specific object types (fibre):

\begin{defi}[$S\C_n$]
  Given a type signature $S$, a natural number $n \in \NN$, and a category $\C$, we 
  define the category $S\C_n$ to be the category with, as objects, diagrams of the form
     \[ n \stackrel{\vectorletter{t}}{\to} T \stackrel{V}{\to} \C \enspace , \]
  written $(V, t_1, \ldots, t_n)$ with $t_i := \vectorletter{t}(i)$, 
where $T$ is a model of $S$,
  the object $V \in \family{\C}{T}$ is a $T$--indexed family of objects of $\C$ and $\vectorletter{t}$ is 
  a vector of elements of $T$ of length $n$.
  A morphism $h$ to another such $(W,\vectorletter{t}) $
  with the same pointing map $\vectorletter{t}$ is given by a morphism $h : V\to W$ in $\family{\C}{T}$
  whose first component is a morphism of $S$-models;
  there are no morphisms between objects with different pointing maps.

We denote by $SU_n:S\C_n \to \Set$ the functor mapping an object $(T,V,\vectorletter{t})$ 
to the underlying set $T$.

  We have a forgetful functor $S\C_n \to \T\C_n$ which forgets the 
 structure of model.
  On the other hand, any model $T$ of $S$ in a set $T$ gives rise to a functor 
  $\family{\C}{T}_n \to S\C_n$, which ``attaches'' the structure of model.
\end{defi}

The meaning of a term $s\in S(n)$ as a natural transformation \[s: 1 \nt SU_n : S\C_n \to \Set\]
is now given by recursion on the structure of $s$:

\begin{defi}[Canonical natural transformation]\label{def:nat_trans_type_indicator}\label{def:canonical_nat_trans}
  Let $s\in S(n)$ be a type of degree $n$. 
  Then $s$ denotes a natural transformation 
       \[s:1 \nt SU_n : S\C_n \to \Set \enspace \]
  defined recursively on the structure of $s$ as follows: for $s = \alpha (a_1,\ldots,a_k)$
  the image of a constructor $\alpha \in S$ we set
  \[s(T,V,\vectorletter{t}) = \alpha (a_1(T,V,\vectorletter{t}),\ldots,a_k(T,V,\vectorletter{t})) \]
  and for $s = m$ with $1\leq m\leq n$  we define
  \[s(T,V,\vectorletter{t}) = \vectorletter{t}(m) \enspace . \]
  We call a natural transformation of the form $s\in S(n)$ \emph{canonical}.
\end{defi}

\begin{defi}[Elementary half-arity]
  We restrict our attention to \emph{elementary} half--arities, 
  which we define 
  analogously to \cite{ahrens_ext_init} as constructed using derivations and products,
  starting from the fibres of the tautological module and the constant singleton module.
  The following clauses define the inductive set of elementary half-arities:
  \begin{itemize}
   \item The constant functor $* : R \mapsto 1$ is an elementary half--arity.

   \item Given any canonical natural transformation $\tau : 1 \nt S U_n$ (cf.\ Definition\ \ref{def:canonical_nat_trans}), 
         the point-wise fibre module with respect to $\tau$ (cf.\ Definition\ \ref{def:fibre_rel_mod_II}) of the tautological module 
         $\Theta_n : R\mapsto (R_n, R_n)$ (cf.\ Definition\ \ref{def:taut_mod_pointed}) is an elementary half--arity of degree $n$, 
        \[ \fibre{\Theta_n}{\tau} : \SigRMon{S} \to \LRMod{n}{S}{\Set} \enspace , \quad R\mapsto \fibre{R_n}{\tau} \enspace . \]

   \item Given any elementary half-arity $M : \SigMon{S} \to \LMod{n}{S}{\Set}$ 
         of degree $n$ and a canonical natural transformation $\tau : 1 \nt S U_n$, 
         the point-wise derivation of $M$ with respect to $\tau$ is an elementary half--arity of degree $n$, 
      \[ M^{\tau} : \SigRMon{S} \to \LRMod{n}{S}{\Set} \enspace , \quad R\mapsto \bigl(M(R)\bigr)^{\tau} \enspace . \]
       Here $\bigl(M(R)\bigr)^{\tau}$ really means derivation of the module, i.e., derivation in the second component
         of $M(R)$.

   \item For a half--arity $M$, let $M_i : R \mapsto \pi_i M(R)$ denote the $i$th projection.
     Given two elementary half--arities $M$ and $N$ of degree $n$, which coincide pointwise on the first 
          component, i.e.\ such that $M_1 = N_1$.
       Then their product $M\times N$ is again an elementary half-arity of degree $n$. Here the product is 
        really the pointwise product in the second component, i.e.\ 
            \[ M\times N : R \mapsto \bigl(M_1(R), M_2(R)\times N_2(R)\bigr) \enspace .  \]
  \end{itemize}

\end{defi}
 
\noindent
As explained at the beginning of Section~\ref{sec:term_sigs}, our goal is to consider 
arities such as application, which are families of arities indexed by object types,
as one arity of higher degree.
The following result explains how this grouping (and ungrouping) translates to half--arities:
a half--arity of degree $n$ thus associates, to any relative $S$--monad $P$ over a set of types $T$, 
a \emph{family of $P$--modules} indexed by $T^n$.

\begin{lem}[Module of higher degree corresponds to a family of modules]
  \label{lem:family_of_mods_cong_pointed_mod}
  \label{lem:family_of_mods_cong_pointed_mod_relative}
 Let $\C$ be a category, let $T$ be a set and $R$ be a monad on $\family{\C}{T}$. 
Suppose $n \in \NN$, and let $\D$ be a category. Then
  modules over $R_n$ with codomain $\D$ correspond precisely to families of $R$--modules indexed by $T^n$ 
  with codomain $\D$ by (un)currying.
 More precisely, let $M$ be an $R_n$--module. Given $\mathbf{t}\in T^n$, we define an $R$--module $M_{\mathbf{t}}$ by 
  \[ M_{\mathbf{t}} (c) := M(c,\mathbf{t}) \enspace . \] 
 Module substitution for $M_{\mathbf{t}}$ is given, for $f \in \family{\C}{T}(c,Rd)$, by
 \[ \mkl[M_{\mathbf{t}}]{f} := \mkl[M]{f} \]
  where we use that we also have $f \in \family{\C}{T}_n ((c,\mathbf{t}), (Rd, \mathbf{t}))$ according to Definition\ \ref{def:cat_set_pointed}.
 Going the other way round, given a family $(M_{\mathbf{t}})_{\mathbf{t}\in T^n}$, we define the $R_n$--module $M$ by
  \[ M(c,\mathbf{t}) := M_{\mathbf{t}}(c) \enspace . \]
 Given a morphism $f \in \family{\C}{T}_n ((c,\mathbf{t}), (Rd, \mathbf{t}))$ --- 
 recall that morphisms in $\family{\C}{T}_n$ are only between families with \emph{the same marker} $\mathbf{t}$ ---, 
we also have $ f \in \family{\C}{T}(c,Rd)$  and define 
 \[ \mkl[M]{f} := \mkl[M_{\mathbf{t}}]{f} \enspace .\]
The remark extends to morphisms of modules; indeed, a morphism of modules $\alpha:M\to N$ on categories with pointed index sets
 corresponds to a family of morphisms $(\alpha_{\vectorletter{t}}:M_{\vectorletter{t}}\to N_{\vectorletter{t}})_{\vectorletter{t}\in T^n}$ 
 between the associated families of modules.
\end{lem}

An arity of degree $n\in \mathbb{N}$ for terms over a type signature $S$ is defined to be a pair of functors
from relative $S$--monads to modules in $\LRMod{n}{S}{\PO}$.
The degree $n$ corresponds to the number of object type indices of its associated constructor.
As an example, the arities of $\abs$ and $\app$ of Display \eqref{eq:sig_tlc_higher_order} are of degree $2$.

\begin{defi}[Weighted set]\label{def:weighted_set}
 A weighted set is a set $J$ together with a map $d:J\to\mathbb{N}$.
\end{defi}

\begin{defi}[Term-arity, signature over $S$]
\label{def:term-arity-semantic}
  An \emph{elementary arity $\alpha$ over $S$ of degree $n$} is a pair 
 \[ s = \bigl(\dom(\alpha), \cod(\alpha)\bigr) \] 
  of half--arities over $S$ of degree $n$ such that 
  \begin{itemize}
   \item $\dom(\alpha)$ is elementary and
   \item $\cod(\alpha)$ is of the form $\fibre{\Theta_n}{\tau}$ for some canonical natural 
        transformation $\tau$ as in Definition \ref{def:canonical_nat_trans}.
  \end{itemize}
Any elementary arity is thus \emph{syntactically} of the form given in Definition\ \ref{def:n_comp_classic_arity_syn}. 
We write $\dom(\alpha) \to \cod(\alpha)$ for the arity $\alpha$, and $\dom(\alpha, R) := \dom(\alpha)(R)$
and similar for the codomain and morphisms of relative $S$--monads. 
Given a weighted set $(J,d)$ as in Definition\ \ref{def:weighted_set},
a \emph{term--signature $\Sigma$ over $S$ indexed by $(J,d)$} is a $J$-family 
 $\Sigma$ of elementary arities over $S$, the arity $\Sigma(j)$ being of degree $d(j)$ for any $j\in J$.
\end{defi}

\begin{defi}[1--signature]\label{def:typed_sig_symantically}
  A \emph{1--signature} is a pair $(S,\Sigma)$ consisting of a type signature $S$ and 
  a term--signature $\Sigma$ (indexed by some weighted set) over $S$.
 By definition, the term-arities of $\Sigma$ are elementary.
\end{defi}

\begin{exa}\label{ex:tlc_sig_higher_order}
The terms of the simply typed lambda calculus over the type signature of Example\ \ref{ex:type_sig_SLC} are given by the arities
 \begin{align*}     
   \abs &: \fibre{(\Theta_2)^1}{2} \to \fibre{\Theta_2}{1\TLCar 2} \enspace , \\
   \app &: \fibre{\Theta_2}{1\TLCar 2} \times \fibre{\Theta_2}{1} \to \fibre{\Theta_2}{2} \quad ,
 \end{align*}
  both of which are of degree $2$.
  Contrast that to the signature for the simply--typed lambda calculus we 
  gave in Display\ \eqref{eq:sig_tlc_simple}. 
  The difference is that here all the family of ``application'' arities, indexed by pairs of object types, are grouped together into 
one arity of degree $2$---and similar for abstraction---whereas this is not the case in Display \eqref{eq:sig_tlc_simple}. 

  Those two arities, $\abs$ and $\app$, can in fact be considered over any type signature $S$ with an arrow constructor, 
  in particular over the signature $S_{\PCF}$ (cf.\ Example\ \ref{ex:term_sig_pcf}).
\end{exa}

\begin{exa}[Example \ref{ex:type_PCF} continued] \label{ex:term_sig_pcf}
  We continue considering $\PCF$. The signature $S_{\PCF}$ for its types 
    is given in Example \ref{ex:type_PCF}.
  The term--signature of \PCF~is given in Figure \ref{fig:pcf_sig_ho}. It consists of an arity for abstraction and an arity for application,
each of degree 2, an arity (of degree 1) for the fixed point operator, and 
one arity of degree 0 for each logic and arithmetic constant, some of which we omit.
  
\begin{figure}[bt]
\centering
\fbox{%
 \begin{minipage}{7cm}

\vspace{-1ex}

 \begin{align*}     
       \abs &: \fibre{(\Theta_2)^1}{2} \to \fibre{\Theta_2}{1\PCFar 2} \enspace , \\
        \app &: \fibre{\Theta_2}{1\PCFar 2} \times \fibre{\Theta_2}{1} \to \fibre{\Theta_2}{2} \enspace ,\\
      \PCFFix &: \fibre{\Theta_1}{1\PCFar 1} \to \fibre{\Theta_1}{1} \enspace , \\
      \PCFn{n} &: * \to \fibre{\Theta}{\Nat }  \qquad \text{for $n\in \NN$}\\
      \PCFSucc &: * \to \fibre{\Theta}{\Nat \PCFar \Nat} \\
      \PCFPred &: * \to \fibre{\Theta}{\Nat \PCFar \Nat} \\
      \PCFZerotest &: * \to \fibre{\Theta}{\Nat \PCFar \Bool } \\
       \PCFcond{\Nat} &: * \to \fibre{\Theta}{\Bool \PCFar \Nat \PCFar \Nat \PCFar \Nat} \\
       \PCFTrue, \PCFFalse &: * \to \fibre{\Theta}{\Bool}\\
	{}&\vdots
  \end{align*}

\vspace{0ex}

\end{minipage}
}

 \caption{Term signature of $\PCF$}\label{fig:pcf_sig_ho}
\end{figure}

\end{exa}

\subsection{Models of 1--signatures}\label{sec:rep_of_1--sigs_typed}

Before giving the general definition, we explain what a model of the simply--typed lambda calculus
          looks like according to our definitions:
          
\begin{exa}
           A model of $(S_{\TLC},\Sigma_{\TLC})$ is given by
          \begin{itemize}\setlength{\itemsep}{.3em}
           \item a model of $S_{\TLC}$ in a set $T$;
           \item a reduction monad $R$ on $\TDelta{T}$;
           \item for any pair $(s,t) \in T^2$, two morphisms of modules
               \[ R^s_t \to R_{s\TLCar t} \quad \text{ and } \quad R_{s \TLCar t} \times R_s \to R_t  \enspace . \]
          \end{itemize}
\end{exa}

\noindent
Note that in the above example, we have \enquote{ungrouped} the operations as explained at the beginning of Section\ \ref{sec:comp_types_sem} 
and Lemma\ \ref{lem:family_of_mods_cong_pointed_mod_relative}.

Most of the work for giving a general definition of model of an arity is already done; the specification is contained in the
definition of signature.

\begin{defi}[Model of an arity, of a signature over $S$]
   \label{def:1--rep_typed}
   A \emph{model of an arity $\alpha$ over $S$ in an $S$--monad $R$} is a morphism of relative modules 
    \[ \dom(\alpha,R) \to \cod(\alpha, R) \enspace . \]
  A \emph{model $R$ of a signature $\Sigma$ over $S$} is a given by a relative $S$--monad---called $R$ as well---and 
  a model $\alpha^R$ of each arity $\alpha$ of $\Sigma$ in $R$.
\end{defi}

Models of $(S,\Sigma)$ are the objects of a category $\Rep^\Delta(S,\Sigma)$, whose morphisms are defined as follows:
\begin{defi}[Morphism of models]
  \label{def:rel_mor_of_reps_typed}
  Given models $P$ and $R$ of a typed signature $(S,\Sigma)$, a morphism of models 
  $f : P\to R$ is given by a morphism of relative $S$--monads $f : P \to R$, such that for any arity $\alpha$ of $S$
  the following diagram of module morphisms commutes:
  \[
  \begin{xy}
   \xymatrix{
        **[l]\dom(\alpha,P) \ar[d]_{\dom(\alpha, f)} \ar[r]^{\alpha^P} & **[r]\cod(\alpha,P) \ar[d]^{\cod(\alpha,f)} \\
        **[l]\dom(\alpha,R) \ar[r]_{\alpha^R} & **[r]\cod(\alpha,R) .
    }
  \end{xy}
 \]
\end{defi}

\begin{lem}\label{lem:init_no_eqs_typed} 
  For any typed signature $(S,\Sigma)$, the category of models of $(S,\Sigma)$ has an initial object.
\end{lem}
\begin{proof}
  The initial object is obtained, analogously to the untyped case (cf.\ \cite{ahrens_relmonads}), 
 via an adjunction $\Delta_* \dashv U_*$ between the categories of models of $(S,\Sigma)$
 in relative monads and those in monads as in \cite{ahrens_ext_init}.
  
  In more detail, to any relative $S$--monad $(T,P) \in \SigRMon{S}$ we associate the $S$--monad
  $U (T,P) := (T,UP)$ where $U_{*}P$ is the monad obtained by postcomposing with the forgetful functor $\family{U}{T} : \TP{T} \to \TS{T}$.
  Substitution for $U_{*}P$ is defined, in each fibre, as in Definition\ \ref{def:rmon_delta_endomon}.
  For any arity $s\in \Sigma$ we have that 
       \[U_{*} \dom(s,P) \cong \dom(s,U_{*}P) \enspace , \] 
  and similar for the codomain.
  The postcomposed model morphism $U_*s(P)$ hence models $s$ in $U_{*}P$ in the sense of \cite{ahrens_ext_init}.
  This defines the functor $U_* : \Rep^\Delta(S,\Sigma) \to \Rep(S,\Sigma)$.
  Conversely, to any $S$--monad we can associate a relative $S$--monad by postcomposing with $\family{\Delta}{T} : \TS{T} \to \TP{T}$,
  analogous to the untyped case in \cite{ahrens_relmonads}, yielding $\Delta_* : \Rep(S,\Sigma) \to \Rep^\Delta(S,\Sigma)$.
  In summary, the natural isomorphism
    \[ \varphi_{R,P}:\bigl(\Rep^\Delta(S,\Sigma)\bigr)(\Delta_* R, P) \cong \bigl(\Rep(S,\Sigma)\bigr)(R, U_*P) \]
  is given by postcomposition with the forgetful functor (from left to right) resp.\ with the functor $\Delta$ (from right to left).
\end{proof}

\begin{rem}[on intrinsic vs.\ extrinsic typing]
 Our approach to signatures and their models uses \emph{intrinsic typing} \cite{dep_syn}.
 This term expresses a method of defining exactly the well--typed terms by organizing them into a family of sets parametrized by 
 object types. It is to be contrasted to \emph{ex}trinsic typing, where one starts out with potentially ill-typed \enquote{preterms},
which then need to be filtered via a well-typedness predicate.
Intrinsic typing has two advantages over extrinsic typing:
\begin{itemize}
 \item Intrinsically typed syntax comes with a more useful recursion principle, which allows for the specification of 
translations between languages that are \emph{automatically} compatible with typing.
 \item Intrinsic typing delegates object level typing to the meta language type system, such as the \coq type system.
                   In this way, the meta level type checker (e.g.\ \coq) sorts out ill--typed terms automatically: writing such a term yields a type error on the meta level.
\end{itemize}

\end{rem}

\subsection{Directed equations}

Analogously to the untyped case (cf.\ \cite{ahrens_relmonads}), a directed equation associates, 
to any model of $(S,\Sigma)$ in a relative monad $P$, two parallel morphisms of $P$--modules.
However, as for arities, a directed equation may be, more precisely, a \emph{family of directed equations}, indexed by object types.
As an example, consider the simply--typed lambda calculus from above, which is defined with \emph{typed} abstraction and application.
Similarly, we have a \emph{typed substitution} operation for $\TLC$, which substitutes a term of type $s\in \TLCTYPE$ for a free variable of type $s$
in a term of type $t\in \TLCTYPE$, yielding again a term of type $t$.
For $s,t\in \TLCTYPE$ and $M\in\SLC(V+s)_t$ and $N \in \SLC(V)_s$, $\beta$-reduction is specified by

  \[ (\lambda_{s,t} M)N \leadsto M [* := N] \enspace , \]
where our notation hides the fact that not only abstraction, but also application and substitution are typed operations.
More formally, such a reduction rule might read as a family of inequalities between morphisms of modules
\[  \comp{(\abs_{s,t} \times\id)}{\app_{s,t}} \enspace \leq \enspace \_ [*^s :=_t \_ ] \enspace , \] 
where $s,t\in \TLCTYPE$ range over types of the simply--typed lambda calculus.
Analogously to Section\ \ref{sec:term_sigs}, we want to specify the $\beta$-rule without referring to the set $\TLCTYPE$,
but instead express it for an arbitrary model $R$ of the typed signature $(S_{\TLC},\Sigma_{\TLC})$ 
(cf.\ Examples\ \ref{ex:type_sig_SLC}, \ref{ex:tlc_sig_higher_order}),
as in

\[  \comp{(\abs^R \times \id)}{\app^R} \enspace \leq \enspace \_ [* := \_ ] \enspace , \] 
where both the left and the right side of the inequality are given by suitable $R$--module morphisms of degree 2.
Source and target of a \emph{half--equation} accordingly are given by functors from models of a typed signature $(S,\Sigma)$
to a suitable category of modules.
A half--equation then is a natural transformation between functors into categories of modules:

\begin{defi}[Category of half-equations] \label{def:half_eq_typed}
Let $(S,\Sigma)$ be a signature. An \emph{$(S,\Sigma)$--module} $U$ of degree $n\in \mathbb{N}$
is a functor from the category of models of $(S,\Sigma)$ as defined in Section\ \ref{sec:rep_of_1--sigs_typed} 
to the category $\LRMod{n}{S}{\PO}$ (cf.\ Definition\ \ref{def:lrmod_typed})
commuting with the forgetful functor to the category of relative monads.
We define a morphism of $(S,\Sigma)$--modules to be a natural transformation which
becomes the identity when composed with the forgetful functor. 
We call these morphisms \emph{half--equations} (of degree $n$).
We write $U^R := U(R)$ for the image of the model $R$ under the $S$--module $U$, and similar for 
morphisms.
\end{defi}

\begin{defi}[Substitution as half-equation]\label{def:subst_half_eq_typed}
  Given a relative monad on $\family{\Delta}{T}$, its associated substitution--of--one--variable operation (cf.\ Definition\ \ref{def:hat_P_subst_typed}) 
  yields a family of module morphisms, indexed by pairs $(s,t)\in T$.
  By Lemma\ \ref{lem:family_of_mods_cong_pointed_mod_relative} this family is equivalent to a module morphism of degree 2.
 The assignment
 \[\subst  : R\mapsto \subst^R :  \fibre{{R}_2^1}{2} \times \fibre{{R}_2}{1} \to \fibre{{R}_2}{2} \]
 thus yields a half--equation of degree $2$ over any signature $S$.
 Its domain and codomain are elementary.
\end{defi} 

\begin{exa}[Example \ref{ex:tlc_sig_higher_order} continued]\label{ex:app_circ_half_typed}
  The map
 \[ \comp{(\abs\times\id)}{\app} : R \mapsto \comp{(\abs^R,\id^R)}{\app^R}  : \fibre{{R}_2^1}{2} \times \fibre{{R}_2}{1} \to \fibre{{R}_2}{2} \]
 is a half--equation over the signature of $\SLC$, as well as over the signature of \PCF.
\end{exa}

\begin{defi}
 \label{def:arity_classic_module_typed}
 Any elementary arity of degree $n$,
  \[s = \fibre{\Theta_n^{{{\tau_1}}}}{\sigma_1} \times \ldots\times \fibre{\Theta_n^{{\tau_m}}}{\sigma_m} \to \fibre{\Theta_n}{\sigma} \enspace , \] 
      defines an elementary $S$--module 
    \[\dom(s) : R\mapsto \fibre{R^{\tau_1}_n}{\sigma_1} \times \ldots \times \fibre{R^{\tau_m}_n}{\sigma_m}  \enspace . \] 
\end{defi}

\begin{defi}[Directed equation, rule] \label{def:ineq_typed}
 Given a signature $(S,\Sigma)$, a \emph{directed equation over $(S,\Sigma)$}, or \emph{$(S,\Sigma)$-rule}, of degree $n\in \NN$ is a pair of 
parallel half--equations between $(S,\Sigma)$--modules of degree $n$.
We write $\alpha \leq \gamma$ for the rule $(\alpha, \gamma)$.
\end{defi}

\begin{exa}[$\beta$-reduction]
  For any 1--signature that has an arity for abstraction and an arity for application,
  we specify $\beta$-reduction using the parallel half--equations of Definition\ \ref{def:subst_half_eq_typed} and Example\ \ref{ex:app_circ_half_typed}:
  \[   \comp{(\abs \times \id)}{\app} \leq \subst : \fibre{\Theta_2^1}{2} \times \fibre{\Theta_2}{1} \to \fibre{\Theta_2}{2} \enspace . \]
\end{exa}

\begin{exa}[Fixpoints and arithmetics of \PCF]\label{ex:pcf_ineqs}
The reduction rules of $\PCF$ are specified by the rules---over the 1-signature of $\PCF$ as given in Example \ref{ex:term_sig_pcf}---of Figure \ref{fig:pcf_reductions_ho}. 
\begin{figure}[bt]
\centering
\fbox{%
 \begin{minipage}{12cm}

\vspace{-1ex}

   \begin{align*}
      \PCFFix &\leq \comp{(\id,\PCFFix)}{\app} : \fibre{\Theta_1}{1\PCFar 1} \to \fibre{\Theta_1}{1} \\
      \comp{(\PCFSucc,\PCFn{n})}{\app} &\leq \PCFn{n+1} : * \to \fibre{\Theta}{\Nat}\\
      \comp{(\PCFPred,\PCFn{0})}{\app} &\leq \PCFn{0} : * \to \fibre{\Theta}{\Nat}\\
      \comp{\left(\PCFPred,\comp{({\PCFSucc},{\PCFn{n}})}{\app}\right)}{\app} &\leq \PCFn{n} : * \to \fibre{\Theta}{\Nat}\\
      \comp{(\PCFZerotest,\PCFn{0})}{\app} &\leq \PCFTrue  : * \to \fibre{\Theta}{\Bool}\\
      \comp{\left(\PCFZerotest,\comp{(\PCFSucc,\PCFn{n})}{\app}\right)}{\app} &\leq \PCFFalse  : * \to \fibre{\Theta}{\Bool}\\
                             &\vdots
   \end{align*}

\vspace{0ex}

\end{minipage}
}

 \caption{Reduction Rules of $\PCF$}\label{fig:pcf_reductions_ho}
\end{figure}

\end{exa}

Given a set $A$ of $(S,\Sigma)$-rules, a model of those rules is any model of $(S,\Sigma)$ which satisfies the rules of $A$ 
in the following sense:

\begin{defi}[Models of rules]\label{def:rep_ineq_typed}
  \label{def:2--rep_typed}
 A \emph{model of an $(S,\Sigma)$-rule $\alpha\leq \gamma : U \to V$} (of degree $n$) is any model 
  $R$ over a set of types $T$ of $(S,\Sigma)$ such that 
  $\alpha^R \leq \gamma^R$ pointwise, i.e., if for any pointed context $(X,\vectorletter{t}) \in \TS{T}\times T^n$, 
      any $t\in T$ and any $y\in U^R_{(X, \vectorletter{t})}(t)$, 
    \begin{equation}\alpha^R(y) \enspace \leq \enspace \gamma^R(y) \enspace , \label{eq:comp_sem_rep_ineq}\end{equation}
where we omit the sort argument $t$ as well as the context $(X,\vectorletter{t})$ from $\alpha$ and $\gamma$.
We say that such a model $R$ \emph{satisfies} the rule $\alpha \leq \gamma$.

For a set $A$ of $(S,\Sigma)$-rules, we call \emph{model of $((S,\Sigma),A)$} any model of $(S,\Sigma)$ that
satisfies each rule of $A$.
We define the category of models of the 2-signature $((S,\Sigma), A)$ to be the full subcategory of the category of
models of $S$ whose objects are models of $((S,\Sigma), A)$.

According to Lemma\ \ref{lem:family_of_mods_cong_pointed_mod_relative}, 
 the inequality of Display \eqref{eq:comp_sem_rep_ineq} is equivalent to ask whether, for any $\vectorletter{t} \in T^n$, 
    any $t\in T$ and any $y\in U_{\vectorletter{t}}^R(X)(t)$, 
\begin{equation*}\alpha_{\vectorletter{t}}^R(y) \enspace \leq \enspace \gamma_{\vectorletter{t}}^R(y) \enspace .
\end{equation*}
\end{defi}

\subsection{Initiality for 2--Signatures}

We are ready to state and prove an initiality result for typed signatures with directed equations:

\begin{thm}\label{thm:init_w_ineq_typed}
 For any set of elementary $(S,\Sigma)$-rules $A$, the category of models of $((S,\Sigma),A)$ has an initial object.
\end{thm}

\begin{proof}
 
The basic ingredients for building the initial model 
  are given by the initial model 
$(\init{S},\init{\Sigma})$---or just $\init{\Sigma}$ for short---in the category $\Rep(S,\Sigma)$ of models in monads on set families 
 \cite{ahrens_ext_init}.
 Equivalently, the ingredients come from the initial object 
$(\init{S},\Delta_*\init{\Sigma})$---or just $\Delta_*\init{\Sigma}$ for short---of models
 \emph{without rules} in the category $\Rep^{\Delta}(S,\Sigma)$ (cf.\ Lemma \ref{lem:init_no_eqs_typed}).
    We call $\init{\Sigma}$ resp.\ $\Delta_*\init{\Sigma}$  the monad resp.\ relative monad underlying the initial model.

  The proof consists of several steps: at first, we define a preorder $\leq_A$ on the terms of $\init{\Sigma}$, induced by the set $A$ of rules.
   Afterwards we show that the data of the model $\init{\Sigma}$---substitution, model morphisms etc.---is 
  compatible with the preorder $\leq_A$ in a suitable sense. This yields a model $\init{\Sigma}_A$ of $(\Sigma,A)$.
  Finally we show that $\init{\Sigma}_A$ is the initial such model.

\begin{description}
 \item [The monad underlying the initial model]

    For any context $X\in \TS{\init{S}}$ and $t\in \init{S}$, we equip $\init{\Sigma}X(t)$ with a preorder $A$ by setting (\emph{morally}, cf.\ below), 
      for $x,y\in \init{\Sigma}X(t)$,

    \begin{equation} x \leq_A y \quad :\Leftrightarrow \quad \forall R : \Rep(\Sigma,A), \quad
                      i_R (x) \leq_R i_R (y) \enspace ,
      \label{eq:order_typed}
    \end{equation}
where $i_R : \Delta_*\init{\Sigma} \to R$ is the initial morphism in the category of models of $(S,\Sigma)$, 
 cf.\ Lemma \ref{lem:init_no_eqs_typed}. 
Note that the above definition in Display \eqref{eq:order_typed} is ill--typed:
 we have $x\in \init{\Sigma}X(t)$, which cannot be applied to (a fibre of) $i_R(X) : \retyping{g}(\init{\Sigma}X) \to R(\retyping{g}X)$.
We denote by $\varphi = \varphi_R$ the natural isomorphism induced by the adjunction 
of Definition\ \ref{def:retyping_adjunction_kan} obtained by retyping along the initial morphism of types $g:\init{S}\to T = T_R$ 
towards the set $T$ of ``types'' of $R$,
\[ \varphi_{X,Y} : \family{\PO}{T}\left(\retyping{g}(\init{\Sigma}X), R(\retyping{g}X)\right) \cong 
                     \family{\PO}{\init{S}}\left(\init{\Sigma}X, \comp{g}{R(\retyping{g}X)}\right) \enspace . \]
Instead of the above definition in Display \eqref{eq:order_typed}, we should really write 
    \begin{equation} x \leq_A y \quad :\Leftrightarrow \quad \forall R : \Rep(\Sigma,A), \quad
                      \left(\varphi(i_{R,X})\right) (x) \leq_R \left(\varphi (i_{R,X})\right) (y) \enspace ,
      \label{eq:order_typed_corrected}
    \end{equation}
where we omit the subscript ``$R$'' from $\varphi$.
We have to show that the map \[ X\mapsto \init{\Sigma}_A X := (\init{\Sigma} X, \leq_A) \] yields a relative monad on $\TDelta{\init{S}}$. 
The missing fact to prove is that the substitution with a morphism 
\[ f\in\TP{\init{S}}(\Delta X, \init{\Sigma}_A Y) \cong \TS{\init{S}}(X,\init{\Sigma} Y)\] 
is compatible with the order $\leq_A$:
given any $f \in \TP{\init{S}}(\Delta X, \init{\Sigma}_A Y)$ we show that 
        \[\sigma^{\init{\Sigma}}(f) \in \TS{\init{S}}(\init{\Sigma} X,\init{\Sigma} Y)\] 
is monotone with respect
to $\leq_A$ and hence (the carrier of) a morphism 
\[\sigma^{\init{\Sigma}_A}(f) \in \TP{\init{S}}(\init{\Sigma}_A X, \init{\Sigma}_A Y) \enspace . \]
We overload the infix symbol $\bind{}{}$ to denote monadic substitution.
Note that this notation now hides an implicit argument giving the sort of the term in which we substitute.
Suppose $x, y \in \init{\Sigma}X(t)$ with $x\leq_A y$, we show 
        \[\bind{x}{f} \enspace \leq_A \enspace \bind{y}{f} \enspace .\] 
Using the definition of $\leq_A$, 
we must show, for a given model $R$ of $(\Sigma,A)$, 
\begin{equation}  \left(\varphi(i_R)\right)(\bind{x}{f}) \enspace \leq_R \enspace \left(\varphi(i_R)\right)(\bind{y}{f}) \enspace .
   \label{eq:proof_of_chap5_1}
 \end{equation}
Let $g$ be the initial morphism of types towards the types of $R$.
Since $i:= i_R$ is a morphism of models and thus in particular a \emph{monad} morphism, 
it is compatible with the substitution of $\init{\Sigma}$ and $R$; we have

\begin{equation}
 \begin{xy}
  \xymatrix @R=4pc @C=5pc{
    \retyping{g}(\init{\Sigma}X) \ar[r]^{\retyping{g}(\sigma(f))} \ar[d]_{i_X} & \retyping{g}(\init{\Sigma}Y) \ar[d]^{i_Y} \\
    R(\retyping{g}X) \ar[r]_{\sigma^R(\comp{\retyping{g}f}{i_Y})}& R (\retyping{g} Y).
}
 \end{xy}
\label{eq:comp_sem_monad_mon_mor_diag}
\end{equation}
By applying the isomorphism $\varphi$ on the diagram of Display \eqref{eq:comp_sem_monad_mon_mor_diag}, we obtain
\begin{align} \comp{\sigma(f)}{\varphi(i_Y)} &= \varphi\left(\comp{\retyping{g}(\sigma(f))}{i_Y}\right) \notag \\
                                              &= \varphi\left(\comp{i_X}{\sigma(\comp{\retyping{g}f}{i_Y})}\right) \notag \\
                                              &= \comp{\varphi(i_X)}{g^*\left(\sigma^R(\comp{\retyping{g}f}{i_Y}) \right)} \enspace . \label{eq:p_chap5_2}
\end{align}
Rewriting the equality of Display \eqref{eq:p_chap5_2} twice in the goal shown in Display \eqref{eq:proof_of_chap5_1} yields the goal
\[ g^*\left(\sigma^R(\comp{\retyping{g}f}{i_Y}) \right)\left((\varphi(i_X))(x) \right) = 
           g^*\left(\sigma^R(\comp{\retyping{g}f}{i_Y}) \right)\left((\varphi(i_X))(y) \right) \enspace ,  \]
which is true since $g^*\left(\sigma^R(\comp{\retyping{g}f}{i_Y}) \right)$ is monotone and
  $(\varphi(i_X))(x) \leq_R (\varphi(i_X))(y)$ by hypothesis.
We hence have defined a monad $\init{\Sigma}_A$ over $\TDelta{\init{S}}$.

It remains to show that this is a \emph{reduction} monad: for $f \leq f'$, we must prove that
$\sigma(f) \leq \sigma(f')$. By Display \eqref{eq:p_chap5_2}, it suffices to show that
\[ {g^*\left(\sigma^R(\comp{\retyping{g}f}{i_Y}) \right)} \leq {g^*\left(\sigma^R(\comp{\retyping{g}f'}{i_Y}) \right)} \]
which follows from the fact that $R$ is a reduction monad.

\end{description}

We need a lemma to proceed:
\begin{lem}\label{lem:useful_lemma_typed}
Given an elementary $S$--module $V : \Rep^{\Delta}(S,\Sigma) \to \LRMod{}{S}{\PO}$ from the category of models of $(S,\Sigma)$  in $S$--monads
  to the large category of modules over $S$--monads and $x,y \in V(\init{\Sigma})(X)(t)$, we have
   \[ x \leq_A y \in V^{\init{\Sigma}}_{X}(t) \quad \Leftrightarrow \quad \forall R : \Rep(S,A), \quad V(i_R)(x)\leq_{V^R_X} V(i_R)(y) \enspace ,\]
 where now and later we omit the arguments $X$ and $(i_R(t))$, e.g., in $V(i_R)(X)(i_R(t))(x)$.
\end{lem}

\begin{proof}[Proof of Lemma \ref{lem:useful_lemma_typed}]

The proof is done by induction on the derivation of ``$V$ elementary''. The only interesting case is where $V = M\times N$ is a product:
    \begin{align*} 
	  (x_1, y_1) \leq (x_2,y_2) &\Leftrightarrow x_1 \leq x_2 \wedge y_1 \leq y_2 \\
                         {}         &\Leftrightarrow  \forall R, M(i_R) (x_1) \leq M(i_R) (x_2) \wedge \forall R, N (i_R) (y_1) \leq N(i_R) (y_2) \\
                         {}         &\Leftrightarrow  \forall R, M(i_R) (x_1) \leq M(i_R) (x_2) \wedge  N (i_R) (y_1) \leq N(i_R) (y_2) \\
                         {}         &\Leftrightarrow  \forall R, V(i_R) (x_1,y_1) \leq V(i_R) (x_2,y_2) \enspace .
      \tag*{\qedhere}
    \end{align*}
\end{proof}

\begin{description}

 \item [Modeling $(S,\Sigma)$ in $\init{\Sigma}_A$]
For the modeling of types there is nothing to do.
Any term arity $s \in \Sigma$ should be modeled by the module morphism $s^{\init{\Sigma}}$, i.e.\ by the model of $s$ in $\init{\Sigma}$. 
We have to show that those models are compatible with the preorder $A$.
Given $x\leq_A y$ in $\dom(s,\init{\Sigma})(X)$, we show (omitting the argument $X$ in $s^{\init{\Sigma}}(X)(x)$)
 \[ s^{\init{\Sigma}} (x) \quad \leq_A \quad s^{\init{\Sigma}}(y) \enspace. \]
By definition, we have to show that, for any model $R$ with initial morphism $i = i_R : \init{\Sigma} \to R$ as before,
\[ \varphi(i_X) (s^{\init{\Sigma}} (x)) \quad \leq_R \quad \varphi(i_X) (s^{\init{\Sigma}}(y)) \enspace. \]
But these two sides are precisely the images of $x$ and $y$ under the upper--right composition of the 
diagram of Definition \ref{def:rel_mor_of_reps_typed} for the morphism of models $i_R$.
By rewriting with this diagram we obtain the goal
\[ s^R \Bigl(\bigl(\dom(s)(i_R)\bigr)(x)\Bigr) \quad \leq_R \quad s^R\Bigl(\bigl(\dom(s) (i_R)\bigr)(y)\Bigr) \enspace. \]
We know that $s^R$ is monotone, thus it is sufficient to show 
\[  \bigl(\dom(s)(i_R)\bigr)(x) \quad \leq_R \quad \bigl(\dom(s) (i_R)\bigr)(y) \enspace. \]
This goal follows from  Lemma\ \ref{lem:useful_lemma_typed} 
(instantiated for the elementary $S$--module $\dom(s)$, cf.\ Definition\ \ref{def:arity_classic_module_typed}) 
and the hypothesis $x \leq_A y$.
We hence have established a model of $(S,\Sigma)$ in $\init{\Sigma}_A$.
From now on we refer to this model by $\init{\Sigma}_A$.

 \item [$\init{\Sigma}_A$ satisfies $A$]

The next step is to show that the model $\init{\Sigma}_A$ satisfies $A$. 
Given a rule
     \[\alpha \leq \gamma : U \to V\] 
of $A$ with an elementary $S$--module $V$, 
   we must show that for any context $X \in \TS{\init{S}}$, any $t\in \init{S}$ and any $x\in U(\init{\Sigma}_A)(X)_t$ in the domain of $\alpha$ we have 
\begin{equation*} \alpha^{\init{\Sigma}_A}(x) \quad \leq_A \quad \gamma^{\init{\Sigma}_A}(x) \enspace , 
\end{equation*}
where here and later we omit the context argument $X$ and the sort argument $t$.
By Lemma \ref{lem:useful_lemma_typed} the goal is equivalent to

\begin{equation}
 \forall R : \Rep(\Sigma,A), \quad V(i_R) (\alpha^{\init{\Sigma}_A}(x)) \quad \leq_{V^R_X} \quad V(i_R) (\gamma^{\init{\Sigma}_A} (x)) \enspace . 
 \label{eq:sigma_a_alt_typed}
\end{equation}
Let $R$ be a model of $(\Sigma,A)$. We continue by proving the statement of Display \ref{eq:sigma_a_alt_typed} for $R$.
    The half-equations $\alpha$ and $\gamma$ are natural transformations from $U$ to $V$;
     since $i_R$ is the carrier of a 
   morphism of $(S,\Sigma)$-models from $\Delta\init{\Sigma}$ to $R$, we can thus  rewrite the goal as
       \[ \alpha^R\bigl(U(i_R)(x)\bigr) \quad \leq_{V^R_X} \quad \gamma^R \bigr(U(i_R)(x)\bigr) \enspace , \]
  which is true since $R$ satisfies $A$.

\item [Initiality of $\init{\Sigma}_A$]

Given any model $R$ of $(\Sigma,A)$, the morphism $i_R$ is monotone with respect to the orders on $\init{\Sigma}_A$ and $R$ by construction of $\leq_A$.
It is hence a morphism of models from $\init{\Sigma}_A$ to $R$.
Uniqueness of the morphisms $i_R$ follows from its uniqueness in the category of models of $(S,\Sigma)$, i.e.\ without rules.
Hence $(\init{S},\init{\Sigma}_A)$ is the initial object in the category of models of $((S,\Sigma),A)$.
\qedhere
\end{description}
\end{proof}

\noindent
In Section~\ref{sec:trans_pcf_ulc} we study an example language in detail:
we construct the syntax of $\PCF$ and show that it yields the initial object in a category of models of $\PCF$.
Afterwards, we use the recursion operator obtained from the initiality property to specify
a translation from $\PCF$ to the untyped lambda calculus.
This translation is semantically faithful with respect to the usual 
reduction relation of $\PCF$---generated by the rules of Example \ref{ex:pcf_ineqs}---and $\beta$-reduction of the lambda calculus.

It is typical for a translation from a high-level language to a low-level language, 
one reduction step in the source is translated to several reduction steps in the target.
Closure under transitivity of our reduction relations is hence crucial for this translation to fit into our framework.

\section{A Translation from \texorpdfstring{$\PCF$}{PCF} to \texorpdfstring{$\ULC$}{ULC} via initiality, in \texorpdfstring{$\mathsf{Coq}$}{Coq}} 
\label{sec:trans_pcf_ulc}

In this section we describe the implementation of the category of models of \PCF, equipped with reduction rules, as
well as of its initial object. This yields an instance of Theorem \ref{thm:init_w_ineq_typed}. 

Before going into the specifics of this example, we unfold, in the next remark, the initiality property to 
make explicit the recipe it gives us to specify translations:

\begin{rem}[Iteration principle by initiality]\label{rem:comp_sem_iteration}
 The universal property of the language generated by a 2--signature yields an \emph{iteration principle}
  to define maps---translations---on this language, which are certified to be compatible with substitution and 
  reduction in the source and target languages. How does this iteration principle work? 
More precisely, what data (and proof) needs to be specified in order to define such a translation via initiality from 
a language, say, $(\init{S},\init{\Sigma}_A)$ to another language $(\init{S}',\init{\Sigma}'_{A'})$, generated by signatures $(S,\Sigma,A)$ 
and $(S',\Sigma',A')$, respectively?
The translation is a morphism---an initial one---in the category of models of the signature $(S,\Sigma,A)$ 
of the source language.
It is obtained by equipping the relative monad $\init{\Sigma}'_{A'}$ underlying the target language with a model
of the signature $(S,\Sigma,A)$. In more detail:
\begin{enumerate}
 \item We give a model of the type signature $S$ in the set $\init{S}'$. By initiality of $\init{S}$, this yields a 
       translation $\init{S} \to \init{S}'$ of sorts.
 \item Afterwards, we specify a model of the term signature $\Sigma$ in the monad $\init{\Sigma}'_{A'}$ by 
        defining suitable (families) of morphisms of $\init{\Sigma}'_{A'}$--modules. This yields a model $R$
           of $(S, \Sigma)$ in the monad $\init{\Sigma}'_{A'}$.
       \newcounter{tempcounter} \setcounter{tempcounter}{\value{enumi}}
\end{enumerate}
By initiality, we obtain a morphism $f : (\init{S},\init{\Sigma}) \to R$ of models of $(S,\Sigma)$, that is,
we obtain a translation from $(\init{S},\init{\Sigma})$ to $(\init{S}',\init{\Sigma}')$ as the colax monad morphism underlying 
the morphism $f$.
However, we have not yet ensured that the translation $f$ is compatible with the respective reduction preorders
in the source and target languages. 
\begin{enumerate}\setcounter{enumi}{\value{tempcounter}}
 \item Finally, we verify that the model $R$ of $(S,\Sigma)$ satisfies the rules of $A$, that is, we check
       whether, for each $\alpha \leq \gamma : U \to V \in A$, and for each context $V$, each $t\in \init{S}$ and $x \in U^R_V(t)$,
        \[    \alpha^R (x) \enspace \leq \enspace \gamma^R (x) \enspace . \]
\end{enumerate}
After verifying that $R$ satisfies the rules of $A$, the model $R$ is in fact a model of $(S,\Sigma,A)$.
The initial morphism $f$ thus yields a faithful translation from $(\init{S},\init{\Sigma}_A)$ to $(\init{S}',\init{\Sigma}'_{A'})$.
\end{rem}

To specify a translation from $\PCF$ to the untyped lambda calculus, we hence follow the steps outlined in Remark~\ref{rem:comp_sem_iteration}.
However, for the implementation in \coq of this instance, we make several simplifications compared to the general theorem:

\begin{itemize}
 \item We do not define a notion of 2-signature, but specify directly a \coq type---and a category---of
       models of \PCF with its reduction rules.
 \item We use dependent \coq types to formalize arities of higher degree (cf.\ Definition \ref{def:half_arity_degree_semantic_typed}),
        instead of relying on modules on pointed categories. A model of an arity of degree $n$
        is thus given by a family of module morphisms (of degree zero), indexed $n$ times over the respective object type as 
        described in Lemma \ref{lem:family_of_mods_cong_pointed_mod_relative}.
 \item The relation on the initial object is not defined via the formula of Display \eqref{eq:order_typed}, but directly 
        through an inductive type, cf.\ Code \ref{code:pcf_eval}, and various closures, cf.\ Code \ref{code:pcf_propag}.
\end{itemize}
Note that a translation from $\PCF$ to the untyped lambda calculus was given already in previous work \cite{ahrens_ext_init}.
That translation was between pure syntax, in the sense that no reduction rules were considered on the source and target,
and consequently the translation was not asked to be faithful with respect to any reductions.
Here, we use the same translation, but upgrade it to a translation between languages equipped with their reduction relations.
To make this article self-contained, we repeat many of the code snippets already given in \cite{ahrens_ext_init}.

\subsection{Formalization of categories}

In this formalization, categories are formalized as ``E-categories'';
this means that a category consists of a type, say, $O$ of objects, and a family, say, $A:O \times O \to \fixedtextsf{Setoid}$ of arrows,
and operations of identity and composition satisfying the usual laws.
Here, a setoid is a type equipped with an equivalence relation denoted \lstinline!a == b! in our formalization, and the axioms of 
a category are stated with respect to that equivalence relation, \emph{not} with respect to the identity type of \coq.
Setoids as morphisms of a category have been used by Aczel \cite{aczel_galois} in LEGO (there a setoid is simply called ``set'')
and Huet and Sa\"ibi \cite{concat} in \coq. A careful analysis of E-categories is given in \cite{palmgren-e-cats}.

When specifying an E-category, one needs to specify, in particular, an equivalence relation on their arrows, that takes the role of
equality of arrows. 
In the category of types and functions of types, a suitable notion of equality of functions is pointwise equality.
In the case of a category of functors, a suitable notion of equality on arrows, that is, on natural transformations, is 
pointwise equality of arrows in the target category. 
This extends to categories of modules: a suitable notion of equality of module morphisms is pointwise equality of the underlying natural transformations.

\begin{rem}[E-categories vs.\ univalent categories]
This article was written before univalent foundations became known to the author.
 At the time of preparing the final version of this article, another approach to categories in type theory has been developed \cite{rezk_completion}:
in univalent foundations, various extensionality principles are available as a consequence of the univalence axiom.
This makes it feasible to use of the Martin-Löf identity type for ``equality of arrows'' in a category. 

Going back to the example of categories of functors and natural transformations mentioned above, 
in univalent foundations one can use propositional extensionality and 
function extensionality to reduce Martin-Löf identity of natural transformations to pointwise Martin-Löf identity of the underlying maps.
One thus recovers the above definition of ``sameness'' for natural transformations: pointwise sameness of the underlying functions.

This extends to equality between morphisms of monads, and to equality between morphisms of modules:
instead of \emph{defining} ``sameness'' to mean pointwise equality of the underlying functions as is done in the setoid-based approach,
one can \emph{show}, in univalent foundations, that Martin-Löf identity is equivalent to pointwise Martin-Löf identity.

For our initiality result this means that both approaches are equivalent, at least in theory.
In the formalization it is difficult to judge the difference without actually formalizing both versions.
While formalizers of category theory often speak of ``setoid hell'' when discussing E-categories,
the advantage of having the right notion of ``sameness'' for morphisms of a category \emph{by definition}---instead of
by theorem as in the univalent foundations---might compensate for some of the unwieldiness encountered in the setoid world.
\end{rem}

In the following, we use a few \coq notations to make the code look more like pen-and-paper mathematics.
This table summarizes some of the notation:
\begin{center}
\begin{tabular}{ll}
\coq notation & meaning \\ \hline
\lstinline!f ;; g! & composition $g \circ f$ in a category (note the differing order) \\
\lstinline! f == g!  & setoid equality of morphisms in a category \\
\lstinline!s ~~> t! & object level (PCF) arrow type\\
\lstinline! f @ x!  & object level application (in a model of PCF) \\
\end{tabular}
\end{center}

\subsection{Models of \texorpdfstring{\PCF}{PCF}}

In this section we explain the formalization of models of $\PCF$ with reduction rules 
(cf.\ Figures \ref{fig:pcf_sig_ho} and \ref{fig:pcf_reductions_ho}).
According to Definitions \ref{def:1--rep_typed} and \ref{def:2--rep_typed}, such a model
consists of 
\begin{enumerate}
   \item a model of the types of \PCF~(in a \coq type \lstinline!U!), cf.\ Example \ref{ex:type_PCF},
   \item a reduction monad \lstinline!P! over the functor $\family{\Delta}{U}$ (in the formalization: \lstinline!IDelta U!) and
   \item a model of the arities of \PCF~(cf.\ Example \ref{ex:term_sig_pcf}), i.e., morphisms of $P$--modules with suitable source and target modules such that
   \item the inequalities defining the reduction rules of \PCF~are satisfied. \label{list:last_item}
\end{enumerate}

\noindent
A model of \PCF hence is a \enquote{bundle}, i.e.\ a record type, whose components---or ``fields''---are these \ref{list:last_item} items. %
We first define what a model of the term signature of \PCF~in a monad $P$ is, 
 in the presence of an $S_\PCF$--monad (cf.\ Definition \ref{def:s-rmon}).
Unfolding the definitions, we suppose given a type \lstinline!Sorts!, a relative monad \lstinline!P! over \lstinline!IDelta Sorts!
and three operations on \lstinline!Sorts!: a binary function \lstinline!Arrow! (denoted by an infixed ``\lstinline!~~>!'')
 and two constants \lstinline!Bool! and \lstinline!Nat!.
\begin{lstlisting}
Variable Sorts : Type.
Variable P : RMonad (IDelta Sorts).
Variable Arrow : Sorts -> Sorts -> Sorts.
Variable Bool : Sorts.
Variable Nat : Sorts.
Notation "a ~~> b" := (Arrow a b) (at level 60, right associativity).
\end{lstlisting}
Here and in the sequel, a short arrow \lstinline!->! denotes the type-theoretic function space in \coq.
A long arrow \lstinline!--->! denotes morphisms in a category; the category at hand will be clear from the context.
Equality of morphisms of a category is denoted by \lstinline!==!.
Furthermore, we write \lstinline!<<! to polymorphically denote the relation on a preordered set.

In this context, a model of \PCF~is given by a bunch of module morphisms satisfying some conditions.
We split the definition into smaller pieces, cf.\ Code \ref{code:rpcf_1-sig} to \ref{code:rpcf_arith}.
Note that \lstinline!M[t]! denotes the fibre module of module \lstinline!M! with respect to \lstinline!t!, 
and \lstinline!d M // u! denotes derivation of module \lstinline!M! with respect to \lstinline!u!.
The module denoted by a star \lstinline!*! is the terminal module, which is the constant singleton module.
\begin{form}[1--Signature of $\PCF$]\label{code:rpcf_1-sig}\hfill
\begin{lstlisting}
Class PCF_model_struct := {
  app : forall u v, (P[u ~~> v]) x (P[u]) ---> P[v]
                       where "A @ B" := (app _ _ _ (A,B));
  abs : forall u v, (d P // u)[v] ---> P[u ~~> v];
  rec : forall t, P[t ~~> t] ---> P[t];
  tttt : * ---> P[Bool];
  ffff : * ---> P[Bool];
  nats : forall m:nat, * ---> P[Nat];
  Succ : * ---> P[Nat ~~> Nat];
  Pred : * ---> P[Nat ~~> Nat];
  Zero : * ---> P[Nat ~~> Bool];
  CondN: * ---> P[Bool ~~> Nat ~~> Nat ~~> Nat];
  CondB: * ---> P[Bool ~~> Bool ~~> Bool ~~> Bool];
  bottom: forall t, * ---> P[t];
  ...
\end{lstlisting}
\end{form}

\noindent
These module morphisms are subject to some directed equations, specifying the reduction rules of $\PCF$.
 The $\beta$-rule reads as
\begin{form}[$\beta$-rule for models of $\PCF$] \label{code:rpcf_beta}\hfill
\begin{lstlisting} 
  beta_red : forall r s V y z, abs r s V y @ z << y[*:= z] ;
  ...
\end{lstlisting}
\end{form}

\noindent
where \lstinline!y[*:= z]! is the substitution of the freshest variable (cf.\ Definition \ref{def:hat_P_subst_typed}) 
as a special case of simultaneous monadic substitution. 
The rule for the fixed point operator says that $\mathbf{Y}(f)$ reduces to $f\left(\mathbf{Y}(f)\right)$:
\begin{form}[Rule for fixpoint operator]\label{code:rpcf_rec}\hfill
\begin{lstlisting}
Rec_A: forall V t g, rec t V g << g @ rec _ _ g 
  ...
\end{lstlisting}
\end{form}

\noindent
The other rules concern the arithmetic and logical constants of \PCF.
Firstly, we have that the conditionals reduce according to the truth value
they are applied to:
\begin{form}[Logical rules of $\PCF$ models]\label{code:rpcf_cond}\hfill
\begin{lstlisting}
  CondN_t: forall V n m, CondN V tt @ tttt _ tt @ n @ m << n ;
  CondN_f: forall V n m, CondN V tt @ ffff _ tt @ n @ m << m ;
  CondB_t: forall V n m, CondB V tt @ tttt _ tt @ n @ m << n ;
  CondB_f: forall V n m, CondB V tt @ ffff _ tt @ n @ m << m ;
   ...
\end{lstlisting}
\end{form}

\noindent
Furthermore, we have that $\Succ(n)$ reduces to $n+1$ (which in \coq is written \lstinline!S n!), 
reduction of the $\zeroqu$ predicate
 according to whether its argument is zero or not, and that the predecessor  is 
 post--inverse to the successor function:
\begin{form}[Arithmetic rules of $\PCF$ models] \label{code:rpcf_arith}\hfill
\begin{lstlisting}
  Succ_red: forall V n, Succ V tt @ nats n _ tt << nats (S n) _ tt ;
  Zero_t: forall V, Zero V tt @ nats 0 _ tt << tttt _ tt ;
  Zero_f: forall V n, Zero V tt @ nats (S n) _ tt << ffff _ tt ;
  Pred_Succ: forall V n, Pred V tt @ (Succ V tt @ nats n _ tt) << nats n _ tt;
  Pred_Z: forall V, Pred V tt @ nats 0 _ tt << nats 0 _ tt  }.
\end{lstlisting}
\end{form}

\noindent
After abstracting over the section variables we package all of this into a record type:

\begin{lstlisting}
Record PCF_model := {
  Sorts : Type;
  Arrow : Sorts -> Sorts -> Sorts;
  Bool : Sorts ;
  Nat : Sorts ;
  pcf_monad :> RMonad (IDelta Sorts);
  pcf_model_struct :> PCF_model_struct pcf_monad  Arrow Bool Nat }.
Notation "a ~~> b" := (Arrow a b) (at level 60, right associativity).
\end{lstlisting}
The type \lstinline!PCF_model! constitutes the type of objects of the category of models of $\PCF$ with reduction rules.

\subsection{Morphisms of models}

A morphism of models (cf.\ Definition\ \ref{def:rel_mor_of_reps_typed}) is built from a morphism $g$ of \emph{type} models
and a colax monad morphism over the retyping functor associated to the map $g$. 
In the particular case of $\PCF$, a morphism of models from $P$ to $R$ consists of a 
morphism of models of the types of \PCF (with underlying map \lstinline!Sorts_map!) and
a colax morphism of relative monads which makes commute the diagrams of the form given in Definition \ref{def:rel_mor_of_reps_typed}.
We first define the diagrams we expect to commute, before packaging everything into a record type of morphisms.
The context is given by the following declarations: 
\begin{lstlisting}
Variables P R : PCF_model.
Variable Sorts_map : Sorts P -> Sorts R.
Hypothesis HArrow : forall u v, Sorts_map (u ~~> v) = Sorts_map u ~~> Sorts_map v.
Hypothesis HBool : Sorts_map (Bool _ ) = Bool _ .
Hypothesis HNat : Sorts_map (Nat _ ) = Nat _ .
Variable f : colax_RMonad_Hom P R
    (RETYPE (fun t => Sorts_map t))
    (RETYPE_PRE (fun t => Sorts_map t))
  (RT_NT (fun t => Sorts_map t)).
\end{lstlisting}
Here the colax monad morphism \lstinline!f! corresponds to the last component of what we defined to 
be a colax monad morphism in Definition \ref{def:colax_rel_mon_mor}, see Remark \ref{rem:rel_mon_mor_case}.

We explain the commutative diagrams of Definition \ref{def:rel_mor_of_reps_typed} for some of the arities. 
For the successor arity we ask the following diagram to commute:

\begin{form}[Commutative diagram for successor arity] \label{code:pcf_succ_diag}\hfill
\begin{lstlisting}
Program Definition Succ_hom :=
    Succ ;; f [(Nat ~~> Nat)] ;; Fib_eq_RMod _ _ ;; IsoPF == *--->* ;; f ** Succ.
\end{lstlisting}
\end{form}

\noindent
Here the morphism \lstinline!Succ! refers to the model of the successor arity either of \lstinline!P! (the
first occurrence) or \lstinline!R! (the second occurrence); \coq~is able to figure this out itself.
The domain of the successor is given by the terminal module $*$. Accordingly, we have that $\dom(\SUCC, f)$ is 
the trivial module morphism with domain and codomain given by the terminal module. We denote this module morphism by \lstinline!*--->*!.
The codomain is given as the fibre of $f$ of type $\Nat \PCFar \Nat$.
The two remaining module morphisms are isomorphisms which do not appear in the informal description.
The isomorphism \lstinline!IsoPF! is needed to permute fibre with pullback (cf.\ Lemma \ref{lem:rel_pb_fibre}).
The morphism \lstinline!Fib_eq_RMod M H! takes a module \lstinline!M! and a proof \lstinline!H! of equality of two object types as 
arguments, say, \lstinline!H : u = v!. Its output is an isomorphism \lstinline! M[u] ---> M[v]!. Here the proof is of type
 \begin{lstlisting}
Sorts_map (Nat ~~> Nat) = Sorts_map Nat ~~> Sorts_map Nat 
 \end{lstlisting}
and \coq is able to figure out the proof itself. 
The diagram for application uses the product of module morphisms, denoted by an infixed \lstinline!X!:
\begin{form}[Commutative diagram for application arity]\label{code:pcf_app_diag}\hfill
\begin{lstlisting}
Program Definition app_hom' := forall u v,
    app u v;; f [( _ )] ;; IsoPF ==
    (f [(u ~~> v)] ;; Fib_eq_RMod _ (HArrow _ _);; IsoPF ) X (f [(u)] ;; IsoPF ) ;; 
       IsoXP ;; f ** (app _ _ ).
\end{lstlisting}
\end{form}

\noindent
In addition to the already encountered isomorphism \lstinline!IsoPF! we have to insert an isomorphism \lstinline!IsoXP! 
which permutes pullback and product (cf.\ Lemma \ref{lem:rel_pb_prod}).
As a last example, we present the property for the abstraction:
\begin{form}[Commutative diagram for abstraction arity]\label{code:pcf_abs_diag}\hfill
\begin{lstlisting}
Program Definition abs_hom' := forall u v,
    abs u v ;; f [( _ )] ==
    DerFib_RMod_Hom _ _ _ ;;  IsoPF ;; f ** (abs (_ u) (_ v)) ;; IsoFP ;; 
       Fib_eq_RMod _ (eq_sym (HArrow _ _ )) .
\end{lstlisting}
\end{form}

\noindent
Here the module morphism \lstinline!DerFib_RMod_Hom f u v! corresponds to the morphism 
\[\dom(\Abs(u,v),f) = \fibre{f^u}{v}, \] and 
\lstinline!IsoFP! permutes fibre with pullback, just like its sibling \lstinline!IsoPF!, but the other way round.

We bundle all those properties into a type class:
\begin{lstlisting}
Class PCF_model_Hom_struct := {
  CondB_hom : CondB_hom' ;
  CondN_hom : CondN_hom' ;
  Pred_hom : Pred_hom' ;
  Zero_hom : Zero_hom' ;
  Succ_hom : Succ_hom' ;
  fff_hom : fff_hom' ;
  ttt_hom : ttt_hom' ;
  bottom_hom : bottom_hom' ;
  nats_hom : nats_hom' ;
  app_hom : app_hom' ;
  rec_hom : rec_hom' ;
  abs_hom : abs_hom' }.
\end{lstlisting}

\noindent
Similarly to what we did for models, we abstract over the section variables and define a record type of 
morphisms of models from \lstinline!P! to \lstinline!R! :
\begin{lstlisting}
Record PCF_model_Hom := {
  Sorts_map : Sorts P -> Sorts R ;
  HArrow : forall u v, Sorts_map (u ~~> v) = Sorts_map u ~~> Sorts_map v;
  HNat : Sorts_map (Nat _ ) = Nat R ;
  HBool : Sorts_map (Bool _ ) = Bool R ;
  model_Hom_monad :> colax_RMonad_Hom P R (RT_NT Sorts_map);
  model_colax_Hom_monad_struct :> PCF_model_Hom_struct
                 HArrow HBool HNat model_Hom_monad }.
\end{lstlisting}

\subsection{Equality of morphisms, category of models}

We have already seen how some definitions that are trivial in informal mathematics, turn into something awful 
in intensional type theory. Equality of morphisms of models is another such definition.
Informally, two such morphisms $a, c : P \to R$ of models are equal if 
\begin{enumerate}
 \item their map of object types $f_a$ and $f_c$ (\lstinline!Sorts_map!) are equal and
 \item their underlying colax morphism of monads---also called $a$ and $c$---are equal.
\end{enumerate}
In our formalization, the second condition is not even directly expressible, since these monad morphisms 
  do not have the same type: 
 we have, for a context $V \in \TS{P}$,
\[ a_V : \retyping{f_a}(PV) \to R(\retyping{f_a}V) \]
and 
\[ c_V : \retyping{f_c}(PV) \to R(\retyping{f_c}V)  \enspace . \]
where $\TS{P}$ is a notation for contexts typed over the set of object types the model $P$ comes with, 
formally the type \lstinline!Sorts P!.
We can only compare $a_V$ to $c_V$ by composing each of them with a suitable transport \lstinline!transp! again, yielding morphisms
\[ \comp{a_V}{R(\text{\lstinline!transp!})} : \retyping{f_a}(PV) \to R(\retyping{f_a}V) \to R(\retyping{f_c}V)\]
and
\[ \comp{\text{\lstinline!transp'!}}{c_V}   : \retyping{f_a}(PV) \to \retyping{f_c}(PV) \to R(\retyping{f_c}V)  \enspace . \]
As before, for equal fibres $\fibre{M}{u}$ and $\fibre{M}{t}$ with $u = t$, the carriers of those transports 
\lstinline!transp! and \lstinline!transp'!
are terms of the form \lstinline!eq_rect _ _ _ H!, where \lstinline!H! is a proof term which depends on 
the proof of 
\begin{lstlisting}
forall x : Sorts P, Sorts_map c x = Sorts_map a x
\end{lstlisting}
of the first condition.
Altogether, the definition of equality of morphisms of models is given by the following inductive proposition:
\begin{lstlisting}
Inductive PCF_Hom_eq (P R : PCF_model) : relation (PCF_model_Hom P R) :=
 | eq_model : forall a c : PCF_model_Hom P R, 
            forall H : (forall t, Sorts_map c t = Sorts_map a t),
            (forall V, a V ;; rlift R (Transp H V)
                                    == 
                      Transp_ord H (P V) ;; c V ) -> PCF_Hom_eq a c.
\end{lstlisting}
The formal proof that the relation thus defined is an equivalence is inadequately long when compared to 
its mathematical complexity, due to the transport elimination.

Composition of models is done by composing the underlying maps of sorts, as well as composing
the underlying monad morphisms pointwise. Again, this operation, which is trivial from a mathematical point
of view, yields a difficulty in the formalization, due to the fact that in \coq,
$\retyping{g}(\retyping{f}V)$ is not convertible to $\retyping{(\comp{f}{g})}V$.
More precisely, suppose given two morphisms of models $a : P \to Q$ and $b : Q \to R$,
given by families of morphisms indexed by $V$ resp.\ $W$,
\begin{align*} a_V &: \widetilde{PV}^a \to Q(\widetilde{V}^a)  \quad \text{and} \\
        b_W &: \widetilde{QW}^b \to R(\widetilde{W}^b) \enspace ,
\end{align*}
where we write $\widetilde{V}^a$ for $\retyping{f_a}V$.
The monad morphism underlying the composite morphism of models is given by the following definition:

\[
 \begin{xy}
    \xymatrix{     \widetilde{PV}^{\comp{a}{b}} \ar[rr]^{\comp{a}{b}_V} \ar[d]_{\text{\lstinline!match!}} & &  R(\widetilde{V}^{\comp{a}{b}}) \\
                        PV  \ar[d]_{\text{\lstinline!ctype!}}                     &  &      R\left(\widetilde{\widetilde{V}^a}^b\right)   \ar[u]_{R(\cong)}  \\
                        \widetilde{PV}^a  \ar[r]_{a_V}        & Q(\widetilde{V}^a) \ar[r]_{\text{\lstinline!ctype!}}& \widetilde{Q(\widetilde{V}^a)}^b \ar[u]_{b_{\widetilde{V}^a}}
       }
 \end{xy}
\]
or, in \coq code,
\begin{lstlisting}
Definition comp_model_car : (forall c : ITYPE U,
        RETYPE (fun t => f' (f t)) (P c) --->
     R ((RETYPE (fun t => f' (f t))) c)) :=
  fun (V : ITYPE U) t (y : retype (fun t => f' (f t)) (P V) t) =>
    match y with ctype _ z =>
      lift (M:=R) (double_retype_1 (f:=f) (f':=f') (V:=V)) _
          (b _ _ (ctype (fun t => f' t)
               (a _ _ (ctype (fun t => f t) z ))))
     end.
\end{lstlisting}
where \lstinline!double_retype_1! denotes the isomorphism in the upper right corner.
The proof of the commutative diagrams for the composite monad morphism is lengthy due to the number of arities of the signature of \PCF.
Definition of the identity morphisms is routine, and in the end we define the category of models of \PCF:
\begin{lstlisting}
Program Instance PCF_MODEL :
     Cat_struct (obj := PCF_model) (PCF_model_Hom) := {
  mor_oid P R := eq_model_oid P R ;
  id R := model_id R ;
  comp a b c f g := model_comp f g }.
\end{lstlisting}

\subsection{One particular model}

We define a particular model, which we later prove to be initial.
First of all, the set of object types of \PCF~is given as follows:
\begin{lstlisting}
Inductive Sorts := 
  | Nat : Sorts
  | Bool : Sorts
  | Arrow : Sorts -> Sorts -> Sorts.
\end{lstlisting}

\noindent
For this section we introduce some notations:
\begin{lstlisting}
Notation "'TY'" := Sorts.
Notation "'IT'" := (ITYPE TY).
Notation "a '~>' b" := (PCF.Arrow a b) (at level 69, right associativity).
\end{lstlisting}

\noindent
We specify the set of \PCF~constants through the following inductive type, indexed by
the sorts of \PCF:
\begin{lstlisting}
Inductive Consts : TY -> Type :=
 | Nats : nat -> Consts Nat
 | ttt : Consts Bool
 | fff : Consts Bool
 | succ : Consts (Nat ~> Nat)
 | preds : Consts (Nat ~> Nat)
 | zero : Consts (Nat ~> Bool)
 | condN: Consts (Bool ~> Nat ~> Nat ~> Nat)
 | condB: Consts (Bool ~> Bool ~> Bool ~> Bool).
\end{lstlisting}
The set family of terms of \PCF~is given by an inductive family, parametrized by a context
\lstinline!V! and indexed by object types: 
\begin{lstlisting}
Inductive PCF (V: TY -> Type) : TY -> Type:=
 | Bottom: forall t, PCF V t
 | Const : forall t, Consts t -> PCF V t
 | Var : forall t, V t -> PCF V t
 | App : forall t s, PCF V (s ~> t) -> PCF V s -> PCF V t
 | Lam : forall t s, PCF (opt t V) s -> PCF V (t ~> s)
 | Rec : forall t, PCF V (t ~> t) -> PCF V t.
Notation "a @ b" := (App a b)(at level 43, left associativity).
Notation "M '" := (Const _ M) (at level 15).
\end{lstlisting}
Monadic substitution is defined recursively on terms:
\lstset{mathescape=false}
\begin{lstlisting}
Fixpoint subst (V W: TY -> Type)(f: forall t, V t -> PCF W t)
           (t : TY)(v : PCF V t) : PCF W t :=
    match v with
    | Bottom t => Bottom W t
    | c ' => c '
    | Var t v => f t v
    | u @ v => u >>= f @ v >>= f
    | Lam t s u => Lam (u >>= shift f)
    | Rec t u => Rec (u >>= f)
    end
where "y >>= f" := (@subst _ _ f _ y).
\end{lstlisting}

\noindent
Here, \lstinline!shift f! is the substitution map \lstinline!f! extended to account for 
an extended context under the binder \lstinline!Lam!.
It is equal to the shifted map of Definition \ref{def:rel_module_deriv}.

\lstset{mathescape=true}
Finally, we define a relation on the terms of type \lstinline!PCF! via the inductive definition

\begin{form}[Reduction rules for $\PCF$] \label{code:pcf_eval}\hfill
\begin{lstlisting}
Inductive eval (V : IT): forall t, relation (PCF V t) :=
 | app_abs : forall (s t:TY) (M: PCF (opt s V) t) N, 
               eval (Lam M @ N) (M [*:= N])
 | condN_t: forall n m, eval (condN ' @ ttt ' @ n @ m) n 
 | condN_f: forall n m, eval (condN ' @ fff ' @ n @ m) m 
 | condB_t: forall u v, eval (condB ' @ ttt ' @ u @ v) u 
 | condB_f: forall u v, eval (condB ' @ fff ' @ u @ v) v
 | succ_red: forall n, eval (succ ' @ Nats n ') (Nats (S n) ')
 | zero_t: eval ( zero ' @ Nats 0 ') (ttt ')
 | zero_f: forall n, eval (zero ' @ Nats (S n)') (fff ')
 | pred_Succ: forall n, eval (preds ' @ (succ ' @ Nats n ')) (Nats n ')
 | pred_z: eval (preds ' @ Nats 0 ') (Nats 0 ')
 | rec_a : forall t g, eval (Rec g) (g @ (Rec (t:=t) g)).
\end{lstlisting}
\end{form}
\noindent
which we then propagate into subterms (cf.\ Code \ref{code:pcf_propag}) and close with respect to transitivity and reflexivity:
\begin{form}[Propagation of reductions into subterms]\label{code:pcf_propag}\hfill
 \begin{lstlisting}
Reserved Notation "x :> y" (at level 70).
Variable rel : forall (V:IT) t, relation (PCF V t).
Inductive propag (V: IT) : forall t, relation (PCF V t) :=
| relorig : forall t (v v': PCF V t), rel v v' -> v :> v'
| relApp1: forall s t (M M' : PCF V (s ~> t)) N, M :> M' -> M @ N :> M' @ N
| relApp2: forall s t (M : PCF V (s ~> t)) N N', N :> N' -> M @ N :> M @ N'
| relLam: forall s t (M M':PCF (opt s V) t), M :> M' -> Lam M :> Lam M'
| relRec: forall t (M M' : PCF V (t ~> t)), M :> M' -> Rec M :> Rec M'
where "x :> y" := (@propag _ _ x y).
 \end{lstlisting}
\end{form}

\noindent
The data thus defined constitutes a relative monad \lstinline!PCFEM! on the functor $\TDelta{T_\PCF}$ (\lstinline!IDelta TY!). 
We omit the details.

Now we need to define a suitable morphism (resp.\ family of morphisms) of \lstinline!PCFEM!--modules for any arity
(of higher degree).
Let $\alpha$ be any such arity, for instance the arity $\App$. We need to verify two things:
\begin{enumerate}
 \item we show that the constructor of \lstinline!PCF! which corresponds to $\alpha$ is monotone with respect to 
        the order on \lstinline!PCFEM!. For instance, we show that for any two terms \lstinline!r s:TY! and any 
          \lstinline!V : IDelta TY!, 
         the function 
         \begin{lstlisting}
fun y => App (fst y) (snd y): PCFEM V (r~>s) x PCFEM V r -> PCFEM V s 
         \end{lstlisting}
is monotone.
 \item We show that the monadic substitution defined above distributes over the constructor,
    i.e.\ we prove that the constructor is the carrier of a \emph{module} morphism.
\end{enumerate}

\noindent
All of these are very straightforward proofs, resulting in a model \lstinline!PCFE_model! of \PCF:
\begin{lstlisting}
Program Instance PCFE_model_struct :
       PCF_model_struct PCFEM PCF.arrow PCF.Bool PCF.Nat := {
  app r s := PCFApp r s;
  abs r s := PCFAbs r s;
  rec t := PCFRec t ;
  tttt := PCFconsts ttt ;
  ffff := PCFconsts fff;
  Succ := PCFconsts succ;
  Pred := PCFconsts preds;
  CondN := PCFconsts condN;
  CondB := PCFconsts condB;
  Zero := PCFconsts zero ;
  nats m := PCFconsts (Nats m);
  bottom t := PCFbottom t }.
Definition PCFE_model : PCF_model := Build_PCF_model  PCFE_model_struct.
\end{lstlisting}
Note that in the instance declaration \lstinline!PCFE_model_struct!, the \lstinline!Program! framework proves 
automatically the properties of Code \ref{code:rpcf_beta}, \ref{code:rpcf_rec}, \ref{code:rpcf_cond} and
 \ref{code:rpcf_arith}.

\subsection{Initiality}\label{sec:comp_sem_formal_initiality}

In this section we define a morphism of models from \lstinline!PCFE_model! to any model \lstinline!R : PCF_model!.
At first we need to define a map between the underlying sorts, that is, a map \lstinline!Sorts PCFE_model -> Sorts R!. 
In short, each \PCF~type goes to its model in \lstinline!R!:
\begin{lstlisting}
Fixpoint Init_Sorts_map (t : Sorts PCFE_model) : Sorts R := 
    match t with
    | PCF.Nat => Nat R 
    | PCF.Bool => Bool R
    | u ~> v => (Init_Sorts_map u) ~~> (Init_Sorts_map v)
    end.
\end{lstlisting}

\noindent
The function \lstinline!init! is the carrier of what will later be proved to be the initial morphism to 
the model \lstinline!R!. It maps each constructor of \PCF~recursively to its counterpart in the
model \lstinline!R!:

\begin{lstlisting}
Fixpoint init V t (v : PCF V t) :
    R (retype (fun t0 => Init_Sorts_map t0) V) (Init_Sorts_map t) :=
  match v with
  | Var t v => rweta R _ _ (ctype _ v)
  | u @ v => app _ _ _ (init u, init v)
  | Lam _ _ v => abs _ _ _ (rlift R
         (@der_comm TY (Sorts R) (fun t => Init_Sorts_map t) _ V ) _ (init v))
  | Rec _ v => rec _ _ (init v)
  | Bottom _ => bottom _ _ tt
  | y ' => match y in Consts t1 return
              R (retype (fun t2 => Init_Sorts_map t2) V) (Init_Sorts_map t1) with
                 | Nats m => nats m _ tt
                 | succ => Succ _ tt
                 | condN => CondN _ tt
                 | condB => CondB _ tt
                 | zero => Zero _ tt
                 | ttt => tttt _ tt
                 | fff => ffff _ tt
                 | preds => Pred _ tt
                 end
  end.
\end{lstlisting}

\noindent
We write $i_V$ for \lstinline!init V! and $g$ for \lstinline!Init_Sorts_map!. 
Note that $i_V : \PCF(V)\to g^*\left(R(\retyping{g}V)\right)$ really is \emph{the image of the initial morphism under
the adjunction $\varphi$} of Definition\ \ref{def:retyping_adjunction_kan}.
Intuitively, passing from \lstinline!init V!$= i_V$ to its adjunct $\varphi^{-1}(i_V)$ is done by 
precomposing with pattern matching on the constructor \lstinline!ctype!.
We informally denote $\varphi^{-1}(i_V)$ by $\comp{\text{\lstinline!match!}}{\text{\lstinline!init V!}}$.

The map \lstinline!init! is compatible with renaming and substitution in \lstinline!PCF! and
\lstinline!R!, respectively, in a sense made precise by the following two lemmas.
The first lemma states that, for any morphism $f : V \to W$ in $\TS{T_\PCF}$, the  following
diagram commutes:
\[
 \begin{xy}
  \xymatrix @C=3.5pc{
        **[l]\PCF(V) \ar[d]_{\text{\lstinline!init V!}} \ar[r]^{\PCF(f)}    &  **[r]\PCF(W) \ar[d]^{\text{\lstinline!init W!}} \\
        **[l]g^*R(\retyping{g}V) \ar[r]_{g^*R(g^*f)}& **[r] g^*R(\retyping{g}W).
}
 \end{xy}
\]

\begin{lstlisting}
Lemma init_lift (V : IT) t (y : PCF V t) W (f : V ---> W) : 
   init (y //- f) = rlift R (retype_map f) _ (init y).
\end{lstlisting}

\noindent
The next commutative diagram concerns substitution; for any $f : V \to \PCF(W)$, 
the diagram obtained by applying $\varphi$ to the diagram given in Display \eqref{eq:comp_sem_monad_mon_mor_diag}---i.e.\ the 
diagram corresponding to Display \eqref{eq:p_chap5_2}---commutes:

\[
 \begin{xy}
  \xymatrix @C=8pc {
        **[l]\PCF(V) \ar[d]_{\text{\lstinline!init V!}} \ar[r]^{\kl[\PCF]{f}}    &  **[r]\PCF(W) \ar[d]^{\text{\lstinline!init W!}} \\
        **[l] g^*R(\tilde{V}) \ar[r]_{g^*\kl[R]{\comp{(g^*f)}{\varphi^{-1}(\text{\lstinline!init W!})}}}& **[r] g^*R(\tilde{W}).
}
 \end{xy}
\]
In \coq the lemma \lstinline!init_subst! proves commutativity of this latter diagram:
\begin{lstlisting}
Lemma init_subst V t (y : PCF V t) W (f : IDelta _ V ---> PCFE W):
  init (y >>= f) =
  rkleisli R (SM_ind (fun t v => match v with ctype t p => init (f t p) end)) 
      _ (init y).
\end{lstlisting}

\noindent
This latter lemma establishes almost the commutative diagram for the family $\varphi^{-1}(i_V)$  to constitute a 
(colax) \emph{monad} morphism, which reads as follows:

\begin{equation}\label{eq:init_subst_monadic}
 \begin{xy}
  \xymatrix @C=5.5pc{
       **[l] \retyping{g}\left(\PCF(V)\right) \ar[d]_{\comp{\text{ \lstinline!match!}}{\text{\lstinline!init V! }}} \ar[r]^{\retyping{g}\left(\kl[\PCF]{f}\right)}    
                           & **[r]\retyping{g}\left({\PCF(W)}\right) \ar[d]^{\comp{\text{ \lstinline!match!}}{\text{\lstinline!init W! }}} \\
        **[l]R(\retyping{g}V) \ar[r]_{\kl[R]{\comp{\comp{(\retyping{g}f)}{\text{ \lstinline!match! }}}{\text{\lstinline!init! }}}}& **[r]R(\retyping{g}{W}) .
}
 \end{xy}
\end{equation}

\noindent
Before we can actually build a monad morphism with carrier map 
  $\comp{\text{\lstinline!match!}}{\text{\lstinline!init V!}}$, we need to verify that \lstinline!init!---and thus 
its adjunct---is monotone.
We do this in 3 steps, corresponding to the 3 steps in which we built up the preorder on the terms
of \PCF:
\begin{enumerate}
 \item the map \lstinline!init! is monotone with respect to the relation \lstinline!eval! (cf.\ Code \ref{code:pcf_eval}):
\begin{lstlisting}
Lemma init_eval V t (v v' : PCF V t) : eval v v' -> init v <<< init v'.
\end{lstlisting}
 \item the map \lstinline!init! is monotone with respect to the propagation into subterms of \lstinline!eval!;
\begin{lstlisting}
Lemma init_eval_star V t (y z : PCF V t) : eval_star y z -> init y <<< init z.
\end{lstlisting}
 \item the map \lstinline!init! is monotone with respect to reflexive and transitive closure of above relation.
\begin{lstlisting}
Lemma init_mono c t (y z : PCFE c t) : y <<< z -> init y <<< init z.
\end{lstlisting}
\end{enumerate}

\noindent
We now have all the ingredients to define the initial morphism from \PCF~to \lstinline!R!.
As already indicated by the diagram in Display \eqref{eq:init_subst_monadic}, its carrier 
is not given by just the map \lstinline!init!, since this map does not have the right type:
its domain is given, for any context $V\in \TS{T_\PCF}$, by $\PCF(V)$ and not, as needed, by $\retyping{g}\left(\PCF(V)\right)$.
We thus precompose with pattern matching in order to pass to its adjunct: 
for any context $V$, the carrier of the initial morphism is given by

\begin{lstlisting}
fun t y => match y with
     | ctype _ p => init p
     end
 : retype _ (PCF V) ---> R (retype _ W)
\end{lstlisting}

\noindent
We recall that the constructor \lstinline!ctype! is the carrier of the natural transformation
of the same name of Definition\ \ref{def:retyping_adjunction_kan}, and that precomposing with pattern matching 
corresponds to specifying maps on a coproduct via its universal property.

\noindent
Putting the pieces together, we obtain a morphism of models of \PCF:
\begin{lstlisting}
Definition initR : PCF_model_Hom PCFE_model R :=
        Build_PCF_model_Hom initR_s.
\end{lstlisting}

\noindent
Uniqueness is proved in the following lemma:
\begin{lstlisting}
Lemma initR_unique : forall g : PCFE_model ---> R, g == initR.
\end{lstlisting}
The proof consists of two steps: first, one has to show that the translation of \emph{sorts}
coincide. Since the source of this translation is an inductive type---the initial model of
 the signature of Example \ref{ex:type_PCF}---this proof is done by induction.
Afterwards the translations of terms are proved to be equal. The proof is done by induction on terms of \PCF. 
It makes essentially use of the commutative diagrams (cf.\ Definition \ref{def:rel_mor_of_reps_typed}) which we presented for 
the arities of successor 
(Code \ref{code:pcf_succ_diag}), application (Code \ref{code:pcf_app_diag}) and abstraction (Code \ref{code:pcf_abs_diag}).
Finally we can declare an instance of \lstinline!Initial! for the category \lstinline!PCF_MODEL! of models:
\begin{lstlisting}
Instance PCF_initial : Initial PCF_MODEL := {
  Init := PCFE_model ;
  InitMor R := initR R ;
  InitMorUnique R := @initR_unique R }.
\end{lstlisting}
Checking the axioms used for the proof of initiality (and its dependencies) yields the use of 
non--dependent functional extensionality (applied to the translations of sorts) and 
uniqueness of identity proofs, which in the \coq standard library is implemented as a consequence of 
another, logically equivalent, axiom \lstinline!eq_rect_eq!:
\begin{lstlisting}
Print Assumptions PCF_initial.
Axioms:
CatSem.AXIOMS.functional_extensionality : forall (A B : Type) (f g : A -> B),
                            (forall x : A, f x = g x) -> f = g
Eq_rect_eq.eq_rect_eq : forall (U : Type) (p : U) (Q : U -> Type) 
                          (x : Q p) (h : p = p), x = eq_rect p Q x p h
\end{lstlisting}

\subsection{A model of \texorpdfstring{\PCF}{PCF} in the untyped lambda calculus}

We use the iteration principle explained in Remark \ref{rem:comp_sem_iteration} in order to 
specify a translation from $\PCF$ to the untyped lambda calculus which is compatible with
reduction in the source and target.
According to the principle, it
is sufficient to define a model of $\PCF$ in the relative 
monad of the lambda calculus 
(cf.\ Example \ref{ex:ulcbeta})
and to verify that this model satisfies the $\PCF$ rules,
formalized in the \coq code snippets \ref{code:rpcf_beta}, \ref{code:rpcf_rec}, \ref{code:rpcf_cond} 
and \ref{code:rpcf_arith}.
The first task, specifying a model of the types of $\PCF$, in the singleton set of types of $\LC$,
is trivial. We furthermore specify models of the term arities of $\PCF$, presented in 
Code \ref{code:rpcf_1-sig}, by giving an instance of the corresponding type class.

\begin{lstlisting}
Program Instance PCF_ULC_model_s :
 PCF_model_struct (Sorts:=unit) ULCBETAM (fun _ _ => tt) tt tt := {
  app r s := ulc_app r s;
  abs r s := ulc_abs r s;
  rec t := ulc_rec t ;
  tttt := ulc_ttt ;
  ffff := ulc_fff ;
  nats m := ulc_N m ;
  Succ := ulc_succ ;
  CondB := ulc_condb ;
  CondN := ulc_condn ;
  bottom t := ulc_bottom t ;
  Zero := ulc_zero ;
  Pred := ulc_pred }.
\end{lstlisting}

\noindent
Before taking a closer look at the module morphisms we specify in order to model the arities of $\PCF$,
we note that in the above instance declaration, we have not given the proofs corresponding to code snippets
\ref{code:rpcf_beta} to \ref{code:rpcf_arith}. 
Referring back to Remark \ref{rem:comp_sem_iteration}, we have not completed the third task, the verification that
the given model satisfies the directed equations. 
The \lstinline!Program! feature we use during the above instance declaration is able to detect that the fields
called \lstinline!beta_red!, \lstinline!rec_A!, etc., are missing, and enters into interactive proof mode to allow us
to fill in each of the missing fields.

We now take a look at some of the lambda terms modeling arities of \PCF.
The carrier of the models of \lstinline!ulc_app! is the application of lambda calculus, of course,
 and similar for \lstinline!ulc_abs!. Here the parameters \lstinline!r! and \lstinline!s! vary over
  terms of type \lstinline!unit!, the type of sorts underlying this model.
We use an infixed application and a de Bruijn notation instead of the more abstract notation of nested data types:
\begin{lstlisting}
Notation "a @ b" := (App a b) (at level 42, left associativity).
Notation "'1'" := (Var None) (at level 33).
Notation "'2'" := (Var (Some None)) (at level 24). 
\end{lstlisting}
 
\noindent
The truth values \True~and \False~are modeled by
\begin{lstlisting}
Eval compute in ULC_True. 
    = Abs (Abs 2)
Eval compute in ULC_False. 
    = Abs (Abs 1)
\end{lstlisting}

\noindent
Natural numbers are given in Church style, the successor function is given by the term $\lambda nfx. f(n~f~x)$.
The predecessor is modeled by the constant 
\[\lambda nfx.n~(\lambda gh.h (g~f)) (\lambda u.x) (\lambda u.u), \]
and the test for zero is modeled by $\lambda n. n (\lambda x.F) T$, where $F$ and $T$ are the lambda terms modeling \False~and \True, respectively.
\begin{lstlisting}
Eval compute in ULC_Nat 0.
    = Abs (Abs 1)
Eval compute in ULC_Nat 2.
    = Abs (Abs (2 @ (Abs (Abs (2 @ (Abs (Abs 1) @ 2 @ 1))) @ 2 @ 1)))
Eval compute in ULC_succ.
    = Abs (Abs (Abs (2 @ (3 @ 2 @ 1))))
Eval compute in ULC_pred.
    = Abs (Abs (Abs (3 @ Abs (Abs (1 @ (2 @ 4))) @ Abs 2 @ Abs 1)))
Eval compute in ULC_zero.
    = Abs (1 @ Abs (Abs (Abs 1)) @ Abs (Abs 2))
\end{lstlisting}

\noindent
The conditional is modeled by the lambda term $\lambda p a b. p~a~b$:
\begin{lstlisting}
Eval compute in ULC_cond.
    = Abs (Abs (Abs (3 @ 2 @ 1)))
\end{lstlisting}

\noindent
The constant arity $\bot_A$ is modeled by $\Omega$:
\begin{lstlisting}
Eval compute in ULC_omega.
    = Abs (1 @ 1) @ Abs (1 @ 1)
\end{lstlisting}

\noindent
The fixed point operator \PCFFix~(\lstinline!rec!) is modeled by the \emph{Turing} fixed--point combinator, that is, the lambda term
\begin{lstlisting}
Eval compute in ULC_theta.
    = Abs (Abs (1 @ (2 @ 2 @ 1))) @ Abs (Abs (1 @ (2 @ 2 @ 1)))
\end{lstlisting}

\noindent
The reason why we use the Turing operator instead of, say, the combinator $\mathbf{Y}$,
\begin{lstlisting}
Eval compute in ULC_Y.
    = Abs (Abs (2 @ (1 @ 1)) @ Abs (2 @ (1 @ 1)))
\end{lstlisting}
is that the latter does not have a property that is crucial for us: we have
\[  \Theta(f) \rightsquigarrow^* f\left(\Theta(f)\right) \]
but only
\[  \mathbf{Y}(f) \stackrel{*}{\leftrightsquigarrow} f\left(\mathbf{Y}(f)\right) \]
via a common reduct. 
Thus if we would attempt to model the arity \lstinline!rec! by the fixed--point combinator $\mathbf{Y}$,
we would not be able to prove the condition expressed in Code \ref{code:rpcf_rec}.
Our more fine-grained approach to operational semantics---using directed equations rather than equations---thus
has the drawback of ruling out the translation of the fixpoint construct of \PCF to the $\mathbf{Y}$ combinator of 
the untyped lambda calculus.

As a final remark, we emphasize that while reduction is given as a relation in our formalization, 
and as such is not computable, the obtained translation from \PCF~to the untyped lambda calculus
is executable in \coq.
For instance, we can translate the \PCF~term negating boolean terms, 
\newcommand{\condB}{\ensuremath{\fixedtextsf{cond}_\Bool}}
$\lambda x.\condB(x)(\false)(\true)$,
 as follows:
\begin{form}[Computing the translation of boolean negation]\label{code:translation_pcf_ulc_example}\hfill
\begin{lstlisting}
Eval compute in 
  (PCF_ULC_c ((fun t => False)) tt (ctype _        
   (Lam (condB ' @@ x_bool @@ fff ' @@ ttt ')))).
   = Abs (Abs (Abs (Abs (3 @ 2 @ 1))) @ 1 @ Abs (Abs 1) @ Abs (Abs 2))
\end{lstlisting}
\end{form}

\noindent
Here the notation ``\lstinline!@@!'' denotes application of \PCF, and \lstinline!x_bool! is
simply a notation for a de Bruijn variable of type \lstinline!Bool! of the lowest level, i.e.\ a variable
that is bound by the \lstinline!Lam! binder of \PCF~in above term.
The resulting term of the lambda calculus is
$\lambda x.( \lambda y. (\lambda z.\lambda w.(w@z@y)) @ x @ (\lambda a.\lambda b. b) @ (\lambda a.\lambda b. a))$.

\section{Conclusion}\label{sec:conclusion}

We have presented an initiality result for simply-typed languages,
based on relative monads from sets to preordered sets.

Our approach is ``intrinsically typed'', that is, terms are typed from the beginning on---there are no ``preterms'', and no predicate of well-typedness.
This yields a useful recursion principle where preservation of typing under translations is guaranteed by construction.

The modeling of reductions via preorders may be considered too coarse: one term might reduce to another term \emph{in different ways},
but the use of preorders to model reduction
                   does not allow to distinguish two reductions with the same source and target.
  
Instead of considering \emph{preordered} sets (indexed by sets of free variables) as models of a 2--signature, 
it would thus be interesting to consider a structure 
which allows for more fine--grained treatment of reduction, such as 
 graphs or categories.
However, one might argue that categories naturally form a bicategory, not a category.
Models of a 2--signature would then also form a bicategory, and one might need to adapt some concepts used here to 
a bicategorical setting.

\subsubsection*{Acknowledgements}
We are very grateful to the anonymous referee for their careful reading and valuable comments, and to the editor Larry Moss for his support.

\bibliographystyle{alpha}
\bibliography{literature}

\end{document}